\documentclass[a4paper,oneside,english]{article}
\usepackage[T1]{fontenc}
\usepackage[latin9]{inputenc}
\usepackage{babel}
\usepackage{units}
\usepackage{amstext,amsmath}
\usepackage{amsthm}
\usepackage{amssymb}
\usepackage{graphicx}
\usepackage{hyperref}

\makeatletter

\newcommand{\email}[1]{\protect\href{mailto:#1}{#1}}

\theoremstyle{plain}
\newtheorem{thm}{\protect\theoremname}
  \theoremstyle{remark}
  \newtheorem{rem}[thm]{\protect\remarkname}

\usepackage{geometry}
\usepackage{tikz,tikz-3dplot}
\usepackage{standalone}

\usepackage{xcolor}
\definecolor{light-gray}{gray}{0.7}

\usepackage{dsfont}

\usepackage[font=small]{caption}
\usepackage{subfig}

\usepackage{amsopn}

\makeatother

  \providecommand{\remarkname}{Remark}

\begin{document}

\global\long\def\supp{{\rm supp}\,}
\global\long\def\R{\mathbb{R}}
\global\long\def\C{\mathbb{C}}
\global\long\def\N{\mathbb{N}}
\global\long\def\Z{\mathbb{Z}}
\global\long\def\epsilon{\varepsilon}
\global\long\def\x{\mathbf{x}}
\global\long\def\sinc{\mathrm{sinc}}
\global\long\def\F{\mathcal{F}}
\global\long\def\E{\mathbb{E}}
 \global\long\def\Var{\textrm{Var}}
 \global\long\def\Tr{\textrm{Tr}}
\global\long\def\Cov{\textrm{Cov}}

\title{Mathematical Analysis of Ultrafast Ultrasound Imaging}

\author{Giovanni S. Alberti\thanks{Department of Mathematics,
ETH Z\"{u}rich, R\"{a}mistrasse 101, 8092 Z\"{u}rich, Switzerland (\email{giovanni.alberti@math.ethz.ch},
\email{habib.ammari@math.ethz.ch}, \email{francisco.romero@sam.math.ethz.ch}).}
\and Habib Ammari\footnotemark[1]
\and Francisco Romero\footnotemark[1]
 \and Timoth\'ee Wintz\thanks{Department of Mathematics and Applications, \'Ecole Normale Sup\'erieure,
45 Rue d'Ulm, 75005 Paris, France (\email{timothee.wintz@ens.fr}).}
}

\date{September 28, 2016}

\maketitle

\begin{abstract}
This paper provides a mathematical analysis of ultrafast ultrasound
imaging. This newly emerging modality for biomedical imaging uses
plane waves instead of focused waves in order to achieve very high frame rates.
We derive the point spread function of the system in the Born approximation
for wave propagation and study its properties. We consider dynamic
data for blood flow imaging, and introduce a suitable random model
for blood cells. We show that a singular value decomposition method
can successfully remove the clutter signal by using the different spatial coherence of tissue and blood signals, thereby providing high-resolution
images of blood vessels, even in cases when the clutter and blood
speeds are comparable in magnitude. Several numerical simulations
are presented to illustrate and validate the approach.
\end{abstract}

\vspace{2mm}


\section{Introduction}

Conventional ultrasound imaging is performed with focused ultrasonic
waves \cite{szabo_diagnostic_2013,shung2015diagnostic}. This yields
relatively good spatial resolution, but clearly limits the acquisition
time, since the entire specimen has to be scanned. Over the last decade,
ultrafast imaging in biomedical ultrasound has been developed \cite{montaldo-tanter-bercoff-benech-finck-2009,62014-tanter-fink,demene2015spatiotemporal}.
Plane waves are used instead of focused waves, thereby limiting the
resolution but increasing  the frame rate considerably, up to 20,000
frames per second. Ultrafast imaging has been made possible by the
recent technological advances in ultrasonic transducers, but the idea
of ultrafast ultrasonography dates back to 1977 \cite{1977-bruneel}.
The advantages given by the very high frame rate are many, and the
applications of this new modality range from blood flow imaging \cite{2011-bercoff,demene2015spatiotemporal},
deep superresolution vascular imaging \cite{Errico2015} and functional
imaging of the brain \cite{Mace2011,2013-mace} to ultrasound elastography
\cite{Gennisson2013487}. In this paper we focus on blood flow imaging.

A single ultrafast ultrasonic image is obtained as follows \cite{montaldo-tanter-bercoff-benech-finck-2009}.
A pulsed plane wave (focused on the imaging plane -- see Figure~\ref{fig:imaging-system})
insonifies the medium, and the back-scattered echoes are measured
at the receptor array, a linear array of piezoelectric transducers.
These spatio-temporal measurements are then beamformed to obtain a
two-dimensional spatial signal. This is what we call \emph{static
inverse problem}, as it involves only a single wave, and the dynamics
of the medium is not captured. The above procedure yields very low
lateral resolution, i.e.\ in the direction parallel to the wavefront,
because of the absence of focusing. In order to solve this issue,
it was proposed to use multiple waves with different angles: these
improve the lateral resolution, but has the drawback of reducing
the frame rate.

For \emph{dynamic imaging}, the above process is repeated many times,
which gives several thousand images per second. In blood flow
imaging, we are interested in locating blood vessels. One of the main
issues lies in the removal of the clutter signal, typically the signal
scattered from tissues, as it introduces major artifacts \cite{2002-bjaerum}.
Ultrafast ultrasonography allows to overcome this issue, thanks to
the very high frame rate. Temporal filters \cite{2011-bercoff,Mace2011,2013-mace},
based on high-pass filtering the data to remove clutter signals, have
shown limited success in cases when the clutter and blood velocities
are close (typically of the order of $\unit[10^{-2}]{m\!\cdot\! s^{-1}}$), or even if the blood velocity is smaller than the clutter velocity. A spatio-temporal method based on the singular value decomposition
(SVD) of the data was proposed in \cite{demene2015spatiotemporal}
to overcome this drawback, by exploiting the different spatial coherence
of clutter and blood scatterers. Spatial coherence is understood as similar movement, in direction and speed, in large parts of the imaged zone. Tissue behaves with higher spatial coherence when compared to the blood flow, since large parts of the medium typically move in the same way, while blood flow is concentrated only in small vessels, which do not share necessarily the same movement direction and speed. This explains why spatial properties are crucial to perform the separation.

In this work, we provide a detailed mathematical
analysis of ultrasound ultrafast imaging. To our knowledge, this
is the first mathematical paper addressing the important challenges
of this emerging and very promising modality. Even though in this work we limit ourselves to formalize the existing methods, the mathematical analysis provided gives important insights, which we expect will lead to improved reconstruction schemes.

The contributions of this paper are twofold. First, we carefully study
the forward and inverse static problems. In particular, we derive
the point spread function (PSF) of the system, in the Born approximation
for ultrasonic wave propagation. We investigate the behavior of the
PSF, and analyze the advantages of angle compounding. In particular,
we study the lateral and vertical resolutions. In addition, this analysis
allows us to fully understand the roles of the key parameters of the
system, such as the directivity of the array and the settings related
to angle compounding.

Second, we consider the dynamic problem. The analysis of the PSF provided
allows to study the doppler effect, describing the dependence on the
direction of the flow. Moreover, we consider a random model for the
movement of blood cells, which allows us to study and justify the SVD
method for the separation of the blood signal from the clutter signal, leading to the reconstruction of the blood vessels' geometry.
The analysis is based on the empirical study of the distribution of
the singular values, which follows from the statistical properties
of the relative data. We provide extensive numerical simulations,
which illustrate and validate this approach.

This paper is structured as follows. In Section~\ref{sect1} we describe the
imaging system and the model for wave propagation. In Section~\ref{sect2} we
discuss the static inverse problem. In particular, we describe the
beamforming process, the PSF and the angle compounding technique.
In Section~\ref{sect3} the dynamic forward problem is considered: we briefly discuss
how the dynamic data are obtained and analyze the doppler effect.
In Section~\ref{sect4} we focus on the source separation to solve the dynamic
inverse problem. We discuss the random model for the refractive index
and the method based on the SVD decomposition of the data. In Section~\ref{sect5}
numerical experiments are provided. Some concluding
remarks and outlooks are presented in the final section.

\section{The Static Forward Problem} \label{sect1}

\begin{figure}
\centering{}\captionsetup[subfigure]{width=175pt}

\subfloat[The real part of the input pulse $f$.\label{fig:pulse}]{\begin{centering}
\includegraphics[bb=0bp 0bp 425bp 318bp,clip,width=0.45\textwidth]{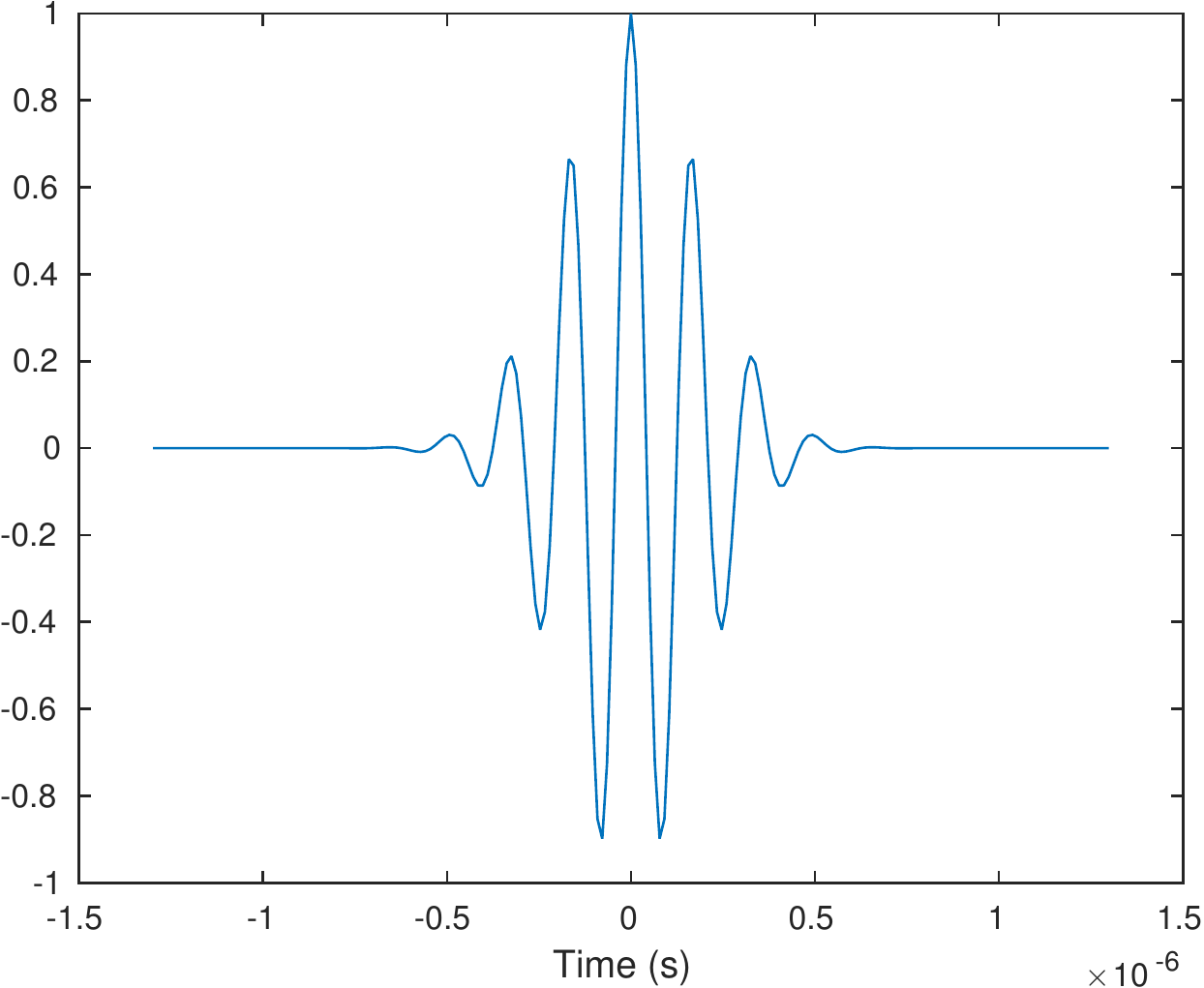}
\par\end{centering}

}\subfloat[\label{fig:imaging-system}The imaging system. The incident wave is
supported near the imaging plane $\{y=0\}$, within the focusing region
bounded by the two curved surfaces.]{\centering{}\includegraphics[width=.52\textwidth]{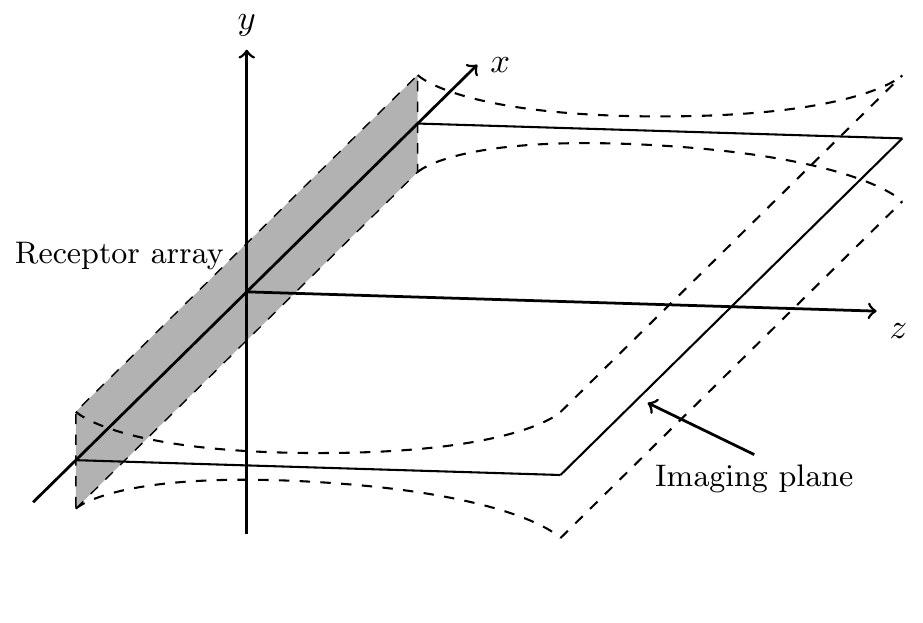}}\caption{The pulse $f$ of the incident wave $u^{i}$ and the focusing region.\label{fig:pulse-system}}
\end{figure}

The imaging system is composed of a medium contained in $\R_{+}^{3}:=\{(x,y,z)\in\R^{3}:z>0\}$
and of a fixed linear array of transducers located on the line $z=0,y=0$.
This linear array of piezoelectric transducers (see \cite[Chapter 7]{szabo_diagnostic_2013})
produces an acoustic illumination that is focused in elevation -- in
the $y$ coordinates, near the plane $y=0$ -- and has the form of
a plane wave in the direction $\mathbf{k}\in S^{1}$ in the $x,z$
coordinates (see Figure~\ref{fig:imaging-system}). Typical sizes for the array length and for the penetration depth are about $\unit[10^{-1}]{m}$.

We make the assumption
that the acoustic incident field $u^{i}$
can be approximated as
\[
u^{i}\left(x,y,z,t\right)=A_{z}\left(y\right)f\left(t-c_{0}^{-1}\mathbf{k}_{}\cdot\left(x,z\right)\right),
\]
where $c_{0}$ is the background speed of sound in the medium. The
function $A_{z}$ describes the beam waist in the elevation direction
at depth $z$ (between $\unit[4\cdot 10^{-3}]{m}$ and $\unit[10^{-2}]{m}$). This is a simplified expression of the true incoming
wave, which is focused by a cylindrical acoustic lens located near
the receptor array (see \cite[Chapters~6 and 7]{szabo_diagnostic_2013}).
The function $f$ is the waveform describing the shape of the input
pulse:
\begin{equation}
f(t)=e^{2\pi i\nu_{0}t}\chi\left(\nu_{0}t\right),\qquad\chi\left(u\right)=e^{-\frac{u^{2}}{\tau^{2}}},\label{eq:f-chi}
\end{equation}
where $\nu_{0}$ is the principal frequency and $\tau$ the width parameter
of the pulse (see Figure~\ref{fig:pulse}). Typically, $\nu_{0}$
will be of the order of $\unit[10^{6}]{s^{-1}}$. More precisely, realistic quantities
are
\begin{equation}\label{eq:quantities}
c_{0}=\unit[1.5\cdot10^{3}]{m\!\cdot\! s^{-1}},\quad\nu_{0}=\unit[6\cdot10^{6}]{s^{-1}},\quad\tau=1.
\end{equation}

Let $c:\R^{3}\rightarrow\mathbb{R}_{+}$ be the speed of sound and
consider the perturbation $n$ given by
\[
n\left(\mathbf{x}\right)=\frac{1}{c^{2}\left(\mathbf{x}\right)}-\frac{1}{c_{0}^{2}}.
\]
We assume that $\mbox{supp}\,n\subseteq\R^3_+$. The acoustic pressure
in the medium satisfies the wave equation
\[
\Delta u\left(\mathbf{x},t\right)-\frac{1}{c^{2}\left(\mathbf{x}\right)}\frac{\partial^{2}}{\partial t^{2}}u\left(\mathbf{x},t\right)=0,\qquad\mathbf{x}\in\mathbb{R}^{3},
\]
with a suitable radiation condition on $u-u^{i}$. Let $G$
denote the Green's function for the acoustic wave equation in $\mathbb{R}^{3}$
\cite{LNM, 2015-watanabe}:
\[
G(\mathbf{x},t,\mathbf{x}',t')=-\frac{(4\pi)^{-1}}{\left|\mathbf{x}-\mathbf{x}'\right|}\delta\left(\left(t-t'\right)-c_{0}^{-1}\left|\mathbf{x}-\mathbf{x}'\right|\right).
\]

In the following, we will assume that the Born approximation holds,
i.e.~we consider only first reflections on scatterers, and neglect
subsequent reflections \cite{LNM, chew1995waves} (in cases when the Born approximation is not valid,  nonlinear methods have to be used). This is a very common approximation in medical imaging, and is justified by the fact that soft biological tissues are almost acoustic homogeneous, due to the high water concentration. In mathematical terms,
it consists in the linearization around the constant sound speed
$c_{0}$. In this case, the scattered wave $u^{s}:=u-u^{i}$ is given
by
\[
u^{s}\left(\mathbf{x},t\right)=\int_{\R}\int_{\mathbb{R}^{3}}n\left(\mathbf{x}'\right)\frac{\partial^{2}u^{i}}{\partial t^{2}}\left(\mathbf{x}',t'\right)G\left(\mathbf{x},t,\mathbf{x}',t'\right)d\mathbf{x}'dt',\qquad\x\in\mathbb{R}^{3},\;t\in\mathbb{R}_{+},
\]
since contributions from $n\partial_{t}^{2}u^{s}$ are negligible.
Therefore, inserting the expressions for the Green's function and
for the incident wave yields

\[
u^{s}\left(\mathbf{x},t\right)=-\int_{\mathbb{R}^{3}}\frac{(4\pi)^{-1}}{\left|\mathbf{x}-\mathbf{x}'\right|}n\left(\mathbf{x}'\right)A_{z'}\left(y'\right)f''\left(t-c_{0}^{-1}\left((x',z')\cdot\mathbf{k}+\left|\mathbf{x}-\mathbf{x}'\right|\right)\right)d\mathbf{x}',
\]
where we set $\x=(x,y,z)$ and $\x'=(x',y',z')$. Since the waist
of the beam in the $y$ direction is small compared to the distance
at which we image the medium, we can make the assumption
\[
\left|\mathbf{x}-\left(x',y',z'\right)\right|\simeq\left|\mathbf{x}-\left(x',0,z'\right)\right|,\qquad \mathbf{x}=\left(x,0,0\right)\in\R^{3},
\]
so that the following expression for $u^{s}$ holds for $\mathbf{x}=\left(x,0,0\right)\in\R^{3}$
and $t>0$:
\[
u^{s}\left(\mathbf{x},t\right)\!=\!\int_{\mathbb{R}^{2}}\frac{-(4\pi)^{-1}}{\left|\mathbf{x}-\left(x',0,z'\right)\right|}f''\left(t-c_{0}^{-1}\left((x',z')\cdot\mathbf{k}+\left|\mathbf{x}-\left(x',0,z'\right)\right|\right)\right)\tilde{n}(x',z')dx'dz',
\]
where $\tilde n$ is given by
\begin{equation}
\tilde{n}(x',z'):=\int_{\mathbb{R}}n\left(\mathbf{x'}\right)A_{z'}\left(y'\right)dy',\qquad\x'=(x',y',z')\in\R^{3}.\label{eq:n-tilde}
\end{equation}

Since our measurements are only two-dimensional (one spatial dimension
given by the linear array and one temporal dimension), we cannot aim
to reconstruct the full three-dimensional refractive index $n$. However,
the above identity provides a natural expression for what can be reconstructed:
the vertical averages $\tilde{n}$ of $n$. 
Since $A_{z}$ is supported near $y=0$, $\tilde{n}$ reflects the
contribution of $n$ only near the imaging plane. In physical terms, $\tilde{n}$ contains all the scatterers in the support of $A_z$; these scatterers are in some sense projected onto $y=0$, the imaging plane. For simplicity,
with an abuse of notation from now on we shall simply denote $\tilde{n}$
by $n$, since the original three-dimensional $n$ will not play any
role, due to the dimensionality restriction discussed above. Moreover,
for the same reasons, all vectors $\x$ and $\x'$ will be two-dimensional,
namely, $\x=(x,z)$ and similarly for $\x'$. In view of these considerations, for $\x=(x,0)\in\R^{2}$ and $t>0$
the scattering wave takes the form 
\begin{equation}
u^{s}\left(\mathbf{x},t\right)=-\int_{\mathbb{R}^{2}}\frac{(4\pi)^{-1}}{\left|\mathbf{x}-\mathbf{x}'\right|}f''\left(t-c_{0}^{-1}\left(\mathbf{x}'\cdot\mathbf{k}+\left|\mathbf{x}-\mathbf{x}'\right|\right)\right)n\left(\mathbf{x}'\right)d\mathbf{x}'.\label{eq:us}
\end{equation}
It is useful to parametrize the direction ${\bf k}\in S^{1}$ of the
incident wave by ${\bf k}={\bf k}_{\theta}=(\sin\theta,\cos\theta)$
for some $\theta\in\R$; in practice, $|\theta|\le0.25$ \cite{montaldo-tanter-bercoff-benech-finck-2009}.

\section{The Static Inverse Problem}  \label{sect2}

The static inverse problem consists in the reconstruction of $n$
(up to a convolution kernel) from the measurements $u^{s}$ at the
receptors, assuming that $n$ does not depend on time. This process
provides a single image, and will be repeated many times in order
to obtain dynamic imaging, as it is discussed in the next sections.

\subsection{Beamforming}

The receptor array is a segment $\Gamma=\left(-A,A\right)\times\left\{ 0\right\} $
for some $A>0$. The travel time from the receptor array to a point
$\mathbf{x}=\left(x,z\right)$ and back to a receptor located in $\mathbf{u}_{0}=\left(u,0\right)$
is given by
\[
\tau_{\mathbf{x}}^{\theta}\left(u\right)=c_{0}^{-1}\left(\mathbf{x}\cdot\mathbf{k}_{\theta}+\left|\mathbf{x}-{\bf u}_{0}\right|\right).
\]
The beamforming process \cite{szabo_diagnostic_2013,montaldo-tanter-bercoff-benech-finck-2009}
consists in averaging the measured signals on $\Gamma$ at $t=\tau_{\mathbf{x}}^{\theta}\left(u\right)$,
which results in the image
\[
s_{\theta}(x,z):=\int_{x-Fz}^{x+Fz}u^{s}\left(\mathbf{u}_{0},\tau_{\mathbf{x}}^{\theta}\left(u\right)\right)du,\qquad\x=(x,z)\in\R_{+}^{2}:=\{(x,z)\in\R^{2}:z>0\}.
\]
The dimensionless aperture parameter $F$ indicates which receptors are chosen to
image the location $\mathbf{x}=\left(x,z\right)$, and depends on
the directivity of the ultrasonic array (in practice,
$0.25\le F\le0.5$ \cite{montaldo-tanter-bercoff-benech-finck-2009}). In general, $F$ depends on the medium roughness and on $\theta$, but  this will not be considered this work. The above identity is the key of the static inverse
problem: from the measurements $u^{s}((u,0),t)$ we reconstruct $s_{\theta}(x,z)$.

\begin{figure}
\captionsetup[subfigure]{width=110pt}\subfloat[The exact PSF given in (\ref{eq:def-g}).]{\includegraphics[clip,width=0.3\textwidth]{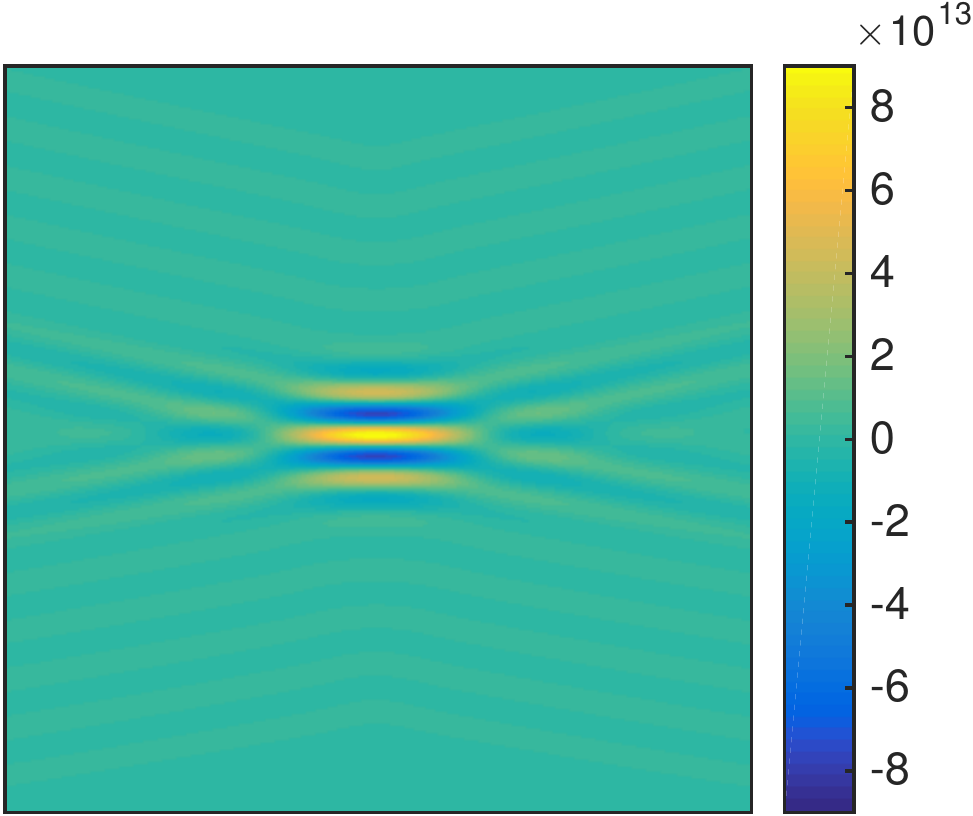}\label{fig:PSF-true}}\hfill{}\subfloat[The approximation of the PSF given in (\ref{eq:g-tilde}).]{\includegraphics[width=0.3\textwidth]{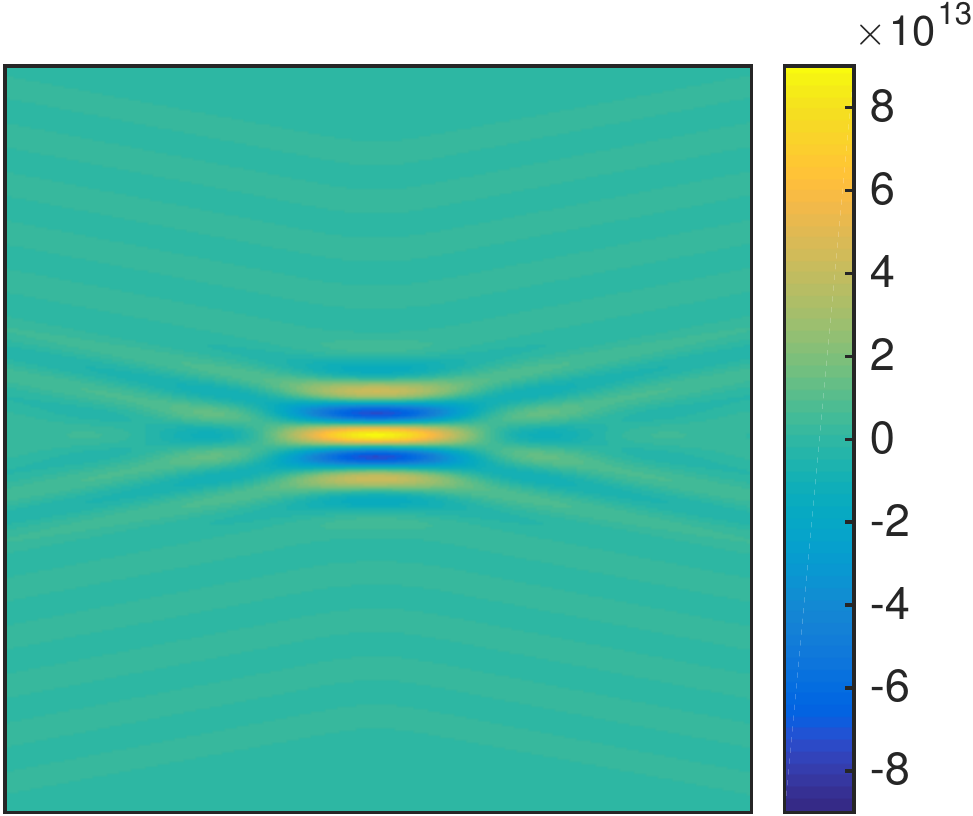}

\label{fig:PSF-app1}}\hfill{}\subfloat[The approximation of the PSF given in (\ref{eq:g-sinc}).]{\includegraphics[width=0.3\textwidth]{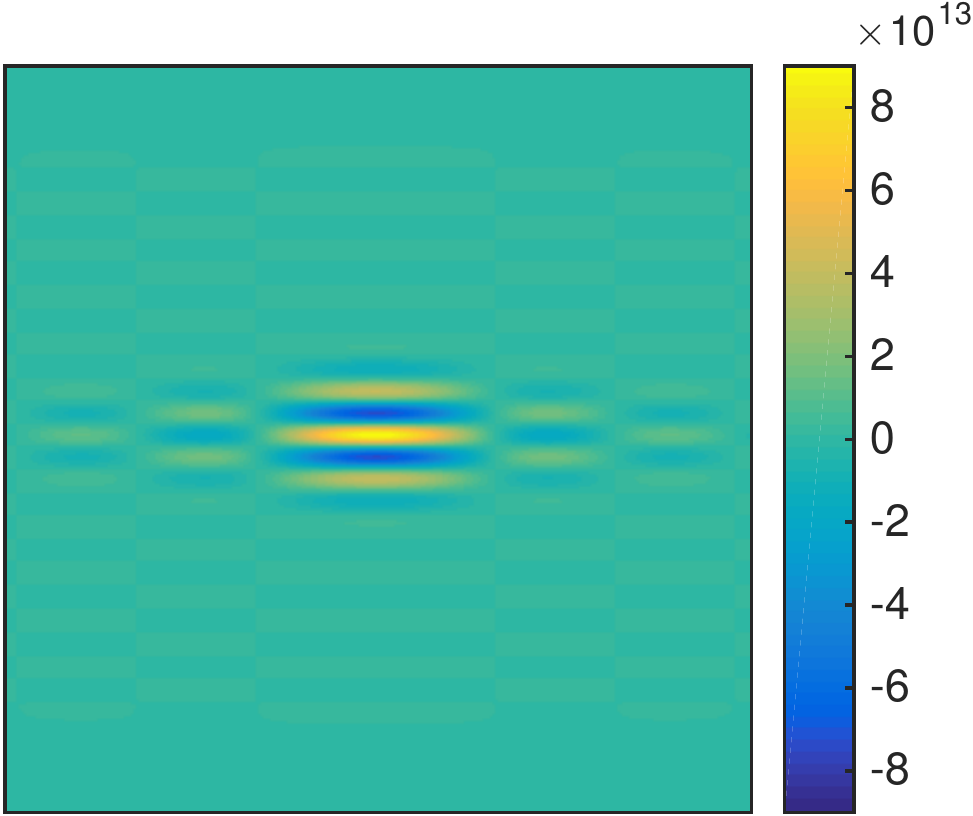}\label{fig:PSF-app2}}\caption{The real part of the point spread function $g_0$  and its approximations
are shown in these figures (with parameters as in \eqref{eq:f-chi} and \eqref{eq:quantities}, and $F=0.4$). The size of the square shown is $\unit[2]{mm}\times\unit[2]{mm}$, and the horizontal and vertical axes are the $x$ and $z$ axes, respectively. The relative error in
the $L^{\infty}$ norm is about 7\% for the approximation shown in
panel (b) and about 9\% for the approximation shown in panel (c).}
\label{fig:PSF}
\end{figure}

We now wish to understand how $s_{\theta}$ is related to $n$. In
order to do so, observe that by (\ref{eq:us}) we may write for $\x\in\R_{+}^{2}$
\begin{equation}\label{eq:s-g-n}
\begin{split}
s_{\theta}(x,z)&=-\int_{\mathbf{x'}\in\mathbb{R}^{2}}\!n\left(\mathbf{x}'\right)\int_{x-Fz}^{x+Fz}\!\frac{(4\pi)^{-1}}{|\x'-{\bf u}_{0}|}f''\left(\tau_{\mathbf{x}}^{\theta}\left(u\right)-\tau_{\mathbf{x}'}^{\theta}\left(u\right)\right)du\,d\mathbf{x}'\\
&=\!\int_{\mathbf{x'}\in\mathbb{R}^{2}}\!g_{\theta}\left(\mathbf{x},\mathbf{x}'\right)n\left(\mathbf{x}'\right)d\mathbf{x}',
\end{split}
\end{equation}
where $g_{\theta}$ is defined as
\begin{equation}
g_{\theta}\left(\mathbf{x},\mathbf{x}'\right)=-\int_{x-Fz}^{x+Fz}\frac{(4\pi)^{-1}}{|\x'-{\bf u}_{0}|}f''\left(\tau_{\mathbf{x}}^{\theta}\left(u\right)-\tau_{\mathbf{x}'}^{\theta}\left(u\right)\right)du,\label{eq:def-g}
\end{equation}
(see Figure~\ref{fig:PSF-true} for an illustration in the case when
$\theta=0$). In other words, the reconstruction $s_{\theta}$ is
the result of an integral operator given by the kernel $g_{\theta}$
applied to the refractive index $n$. Thus, the next step is the study
of the point spread function (PSF) $g_{\theta}\left(\mathbf{x},\mathbf{x}'\right)$,
which should be thought of as the image corresponding to a delta scatterer
in $\x'$.

\subsection{The point spread function}

In its exact form, it does not seem possible to simplify the expression
for $g$ further: we will have to perform some approximations. First,
observe that setting $h_{\x,\x'}^{\theta}(u)=\tau_{\mathbf{x}}^{\theta}\left(u\right)-\tau_{\mathbf{x}'}^{\theta}\left(u\right)$
for  $\x,\x'\in\R_{+}^{2}$ we readily derive
\[
(h_{\x,\x'}^{\theta})'(u)=c_{0}^{-1}(\frac{u-x}{|\x-{\bf u}_{0}|}-\frac{u-x'}{|\x'-{\bf u}_{0}|})\approx c_{0}^{-1}(\frac{u-x}{|\x'-{\bf u}_{0}|}-\frac{u-x'}{|\x'-{\bf u}_{0}|})=c_{0}^{-1}\frac{x'-x}{|\x'-{\bf u}_{0}|},
\]
for $\x$ close to $\x'$ (note that, otherwise, the magnitude of
the PSF would be substantially lower). As a consequence, by (\ref{eq:def-g})
we have
\begin{equation}
\begin{split}g_{\theta}\left(\mathbf{x},\mathbf{x}'\right) & \approx \frac{c_{0}(4\pi)^{-1}}{x-x'}\int_{x-Fz}^{x+Fz}(h_{\x,\x'}^{\theta})'(u)f''\left(h_{\x,\x'}^{\theta}(u)\right)du\\
 & =\frac{c_{0}(4\pi)^{-1}}{x-x'}\left[f'(h_{\x,\x'}^{\theta}(x+Fz))-f'(h_{\x,\x'}^{\theta}(x-Fz))\right].
\end{split}
\label{eq:g-app1}
\end{equation}
In order to simplify this expression even further, let us do a Taylor
expansion of $w_{\pm}^{\theta}(x,z):=h_{\x,\x'}^{\theta}(x\pm Fz)$
with respect to $(x,z)$ around $(x',z')$. Direct calculations show
that
\[
w_{\pm}^{\theta}(x',z')=0,\qquad\nabla w_{\pm}^{\theta}(x',z')=\frac{c_{0}^{-1}}{\sqrt{1+F^{2}}}(\sqrt{1+F^{2}}\sin\theta\mp F,1+\sqrt{1+F^{2}}\cos\theta),
\]
whence
\[
h_{\x,\x'}^{\theta}(x\pm Fz)\!\approx\!\frac{c_{0}^{-1}}{\sqrt{1+F^{2}}}\left(\!(1+\sqrt{1+F^{2}}\cos\theta)(z-z')\!+\!(\sqrt{1+F^{2}}\sin\theta\mp F)(x-x')\!\right)\!.
\]
Substituting this expression into (\ref{eq:g-app1}) yields
\begin{equation}
g_{\theta}(\x,\x')\approx\tilde{g}_{\theta}(\x-\x'),\label{eq:g-conv}
\end{equation}
where
\begin{multline}
\tilde{g}_{\theta}(\x)=\frac{c_{0}}{4\pi x}\left[f'\left(\frac{c_{0}^{-1}}{\sqrt{1+F^{2}}}\left((1+\sqrt{1+F^{2}}\cos\theta)z+(\sqrt{1+F^{2}}\sin\theta-F)x\right)\right)\right.\\
\left.-f'\left(\frac{c_{0}^{-1}}{\sqrt{1+F^{2}}}\left((1+\sqrt{1+F^{2}}\cos\theta)z+(\sqrt{1+F^{2}}\sin\theta+F)x\right)\right)\right],\label{eq:g-tilde}
\end{multline}
(see Figure~\ref{fig:PSF-app1} for an illustration in the case $\theta=0$),
thereby allowing to write the image $s_{\theta}$ given in (\ref{eq:s-g-n})
as a convolution of $\tilde{g}_{\theta}$ and the refractive index
$n$, namely
\[
s_{\theta}(\x)=\int_{\mathbf{x'}\in\mathbb{R}^{2}}\tilde{g}_{\theta}(\x-\x')n\left(\mathbf{x}'\right)d\mathbf{x}'=(\tilde{g}_{\theta}*n)(\x),\qquad\x\in\R_{+}^{2}.
\]

The validity of this approximation, obtained by truncating the Taylor
expansion of $w_{\pm}^{\theta}$ at the first order, is by no means
obvious. Indeed, by construction, the pulse $f(t)$ is highly oscillating ($\nu_0\approx \unit[6\cdot 10^6]{s^{-1}}$),
and therefore even small variations in $t$ may result in substantial
changes in $f(t)$. However, this does not happen, since if $(x,z)$
is not very close to $(x',z')$ then the magnitude of the PSF is very
small, if compared to the maximum value. The verification of this
fact is quite technical, and thus is omitted: the details may be found
in Appendix~\ref{sec:The-justification-of}.
\begin{rem}
From this expression, it is easy to understand the role of the aperture
parameter $F$, which depends on the directivity of the array. Ignoring
the second order effect in $F$ and taking, for simplicity $\theta=0$,
we can further simplify the above expression as
\[
\tilde{g}_{0}(\x)\approx\frac{c_{0}}{4\pi x}\left[f'\left(c_{0}^{-1}\left(2z-Fx\right)\right)-f'\left(c_{0}^{-1}\left(2z+Fx\right)\right)\right].
\]
It is clear that $F$ affects the resolution in the variable $x$:
the higher $F$ is, the higher the resolution is. Moreover, the aperture
parameter affects also the orientation of the diagonal tails in the
PSF. These two phenomena can be clearly seen in Figure~\ref{fig:PSF-F}.
In general, the higher the aperture is the better for the reconstruction:
as expected, the intrinsic properties of the array affects the reconstruction.

\begin{figure}
\captionsetup[subfigure]{width=85pt}\subfloat[The PSF with $F=0.2$.]{\includegraphics[clip,width=0.24\textwidth]{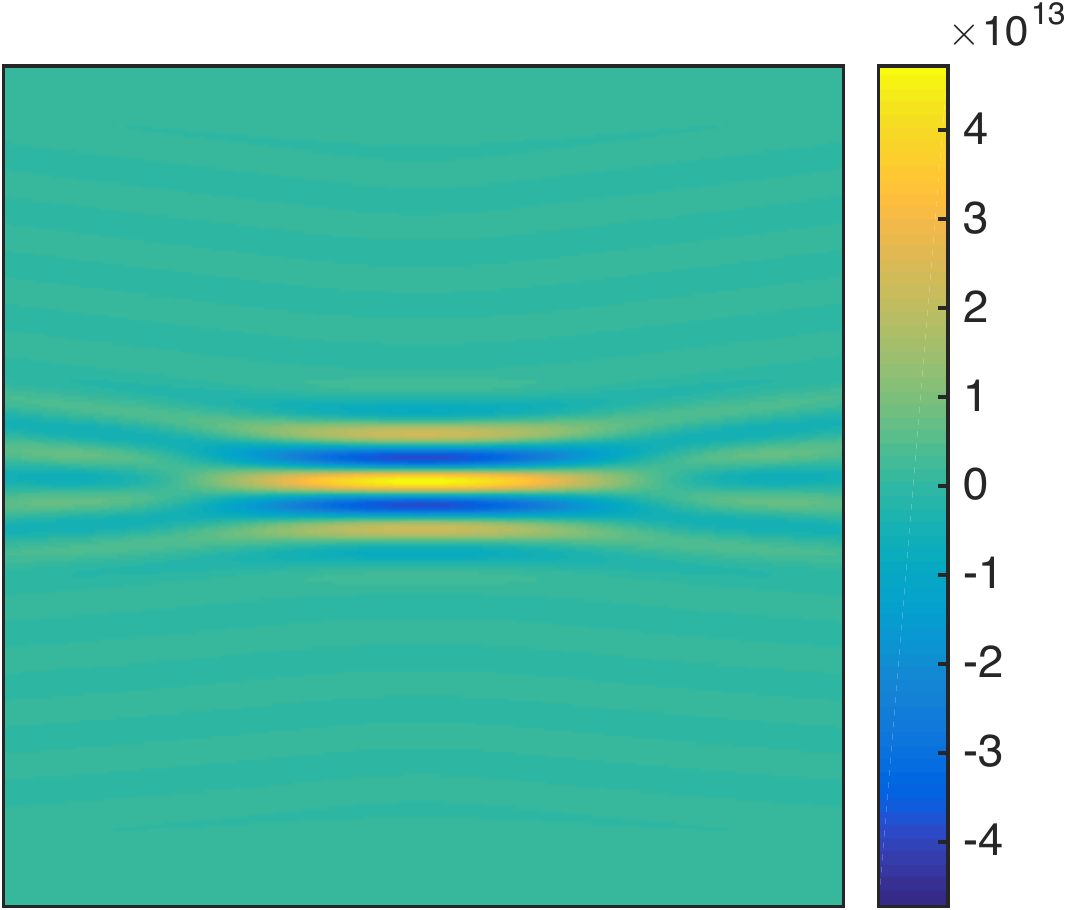}\label{fig:PSF-true-F02}}\hfill{}\subfloat[The PSF with $F=0.3$.]{\includegraphics[width=0.24\textwidth]{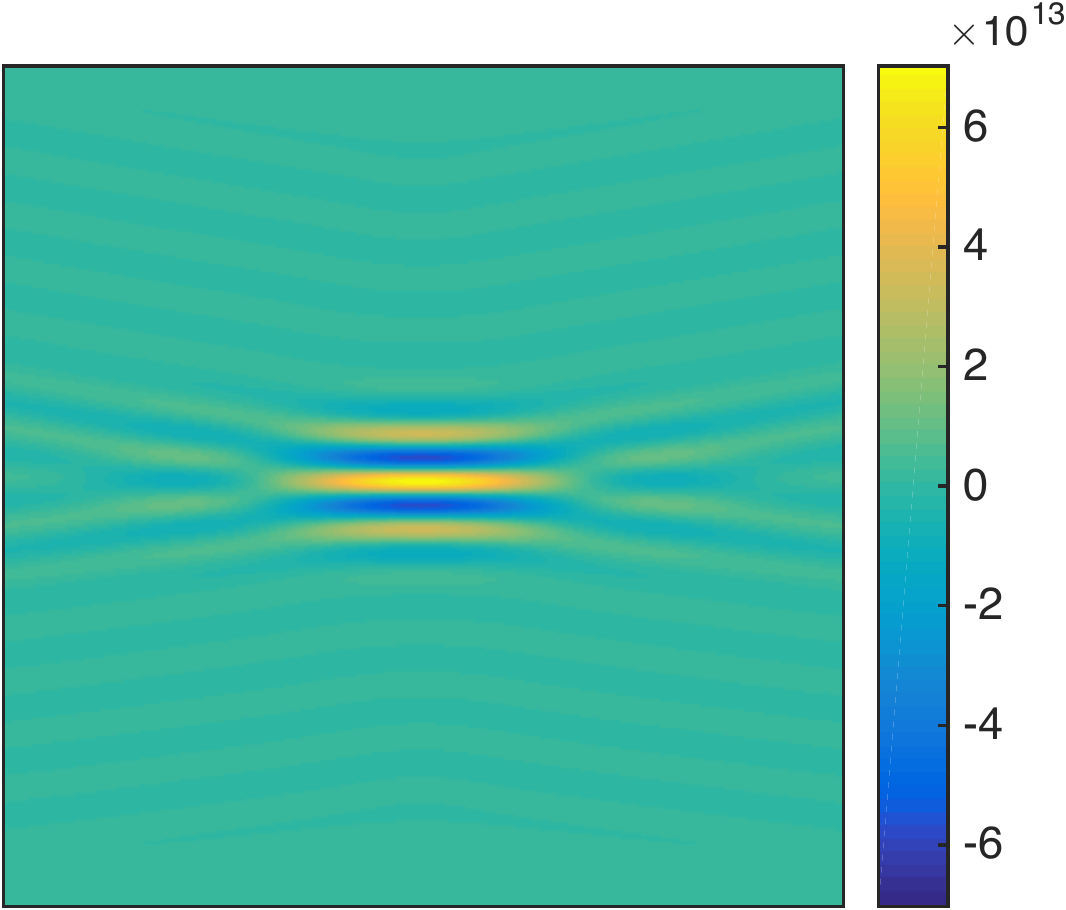}

\label{fig:PSF-true-F03}}\hfill{}\subfloat[The PSF with $F=0.4$.]{\includegraphics[width=0.24\textwidth]{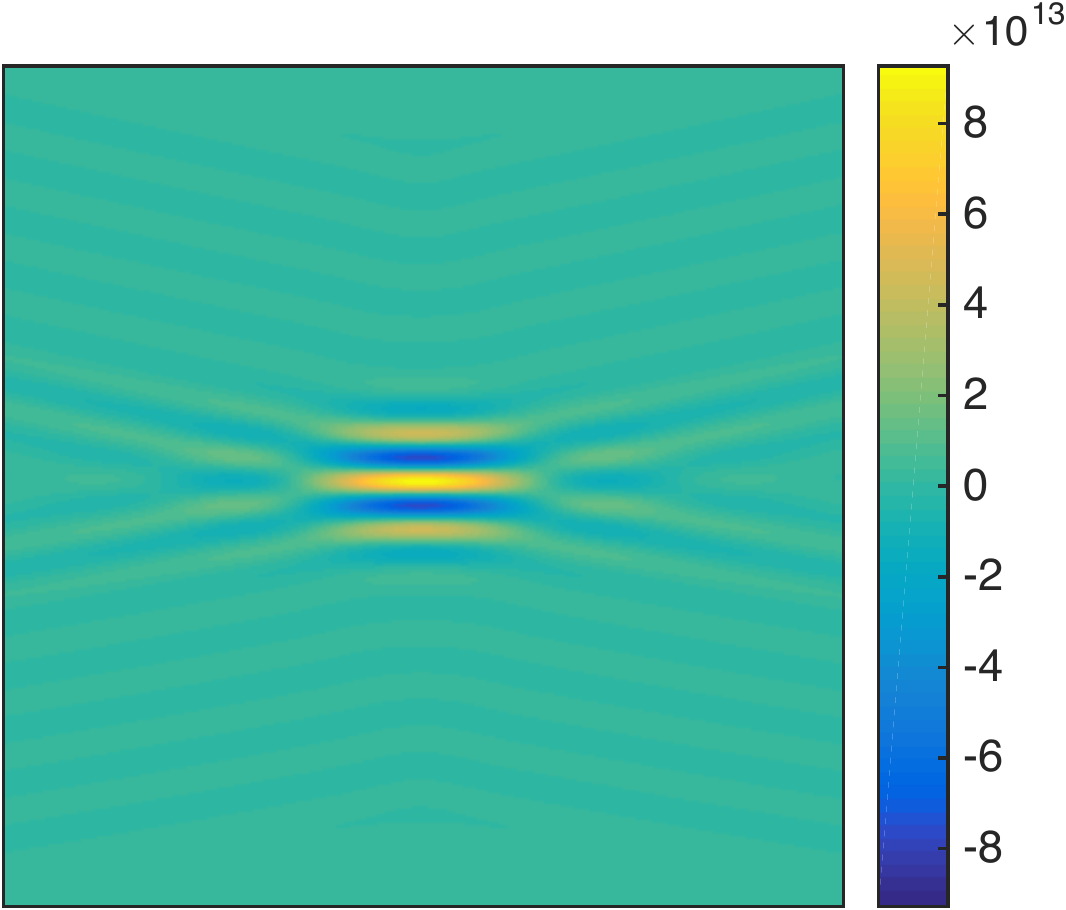}\label{fig:PSF-true-F04}}\hfill{}\subfloat[The PSF with $F=0.5$.]{\includegraphics[width=0.24\textwidth]{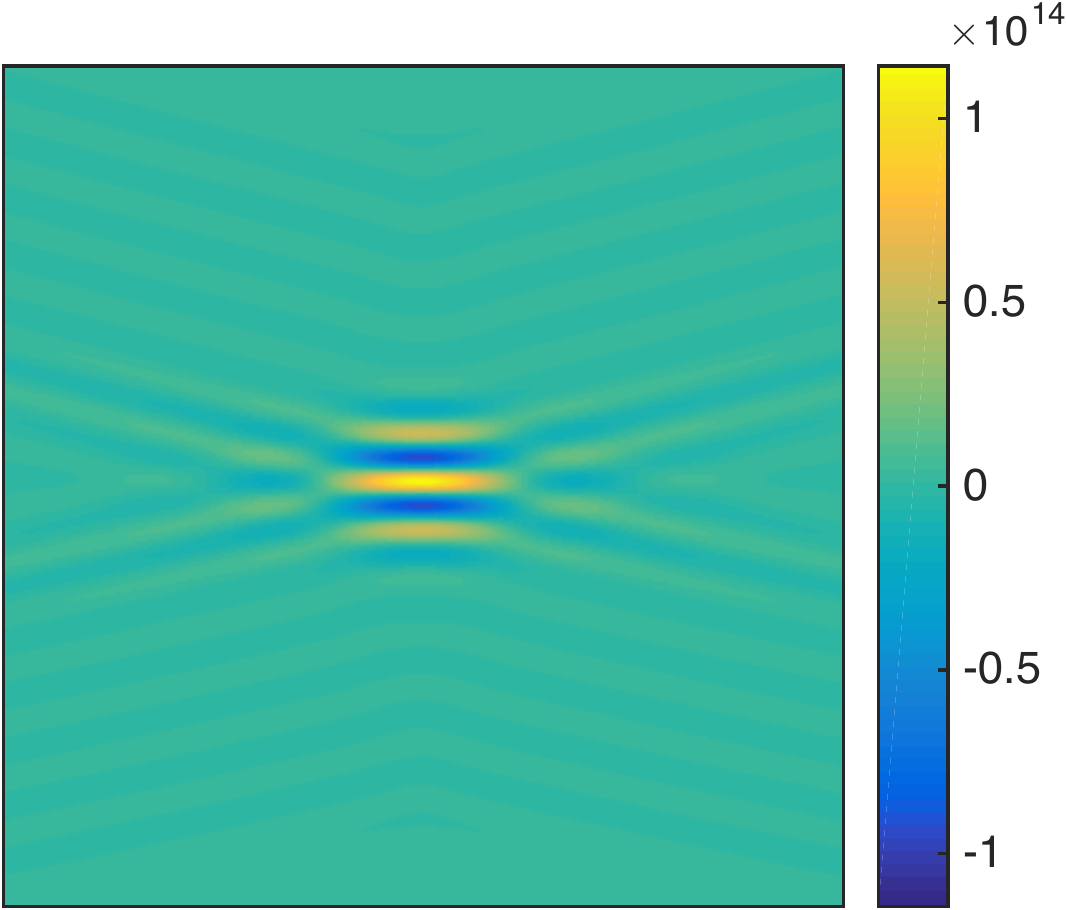}\label{fig:PSF-true-F05}}\caption{The exact PSF with different values of the aperture parameter $F$
(with parameters as in \eqref{eq:f-chi} and \eqref{eq:quantities}, and $\theta=0$). The size of the square shown is $\unit[2]{mm}\times\unit[2]{mm}$, and the horizontal and vertical axes are the $x$ and $z$ axes, respectively.}
\label{fig:PSF-F}
\end{figure}

\end{rem}

\begin{rem}
It is also easy to understand the role of the angle $\theta$. In
view of
\begin{multline*}
\tilde{g}_{\theta}(\x)\approx\frac{c_{0}}{4\pi x}\left[f'\left(c_{0}^{-1}\left((1+\cos\theta)z+(\sin\theta-F)x\right)\right)\right.\\
\left.-f'\left(c_{0}^{-1}\left((1+\cos\theta)z+(\sin\theta+F)x\right)\right)\right],
\end{multline*}
an angle $\theta\neq0$ substantially gives a rotation
of the PSF; see Figure~\ref{fig:PSF-T}.

\begin{figure}
\captionsetup[subfigure]{width=85pt}\subfloat[The PSF with $\theta=0$.]{\includegraphics[clip,width=0.24\textwidth]{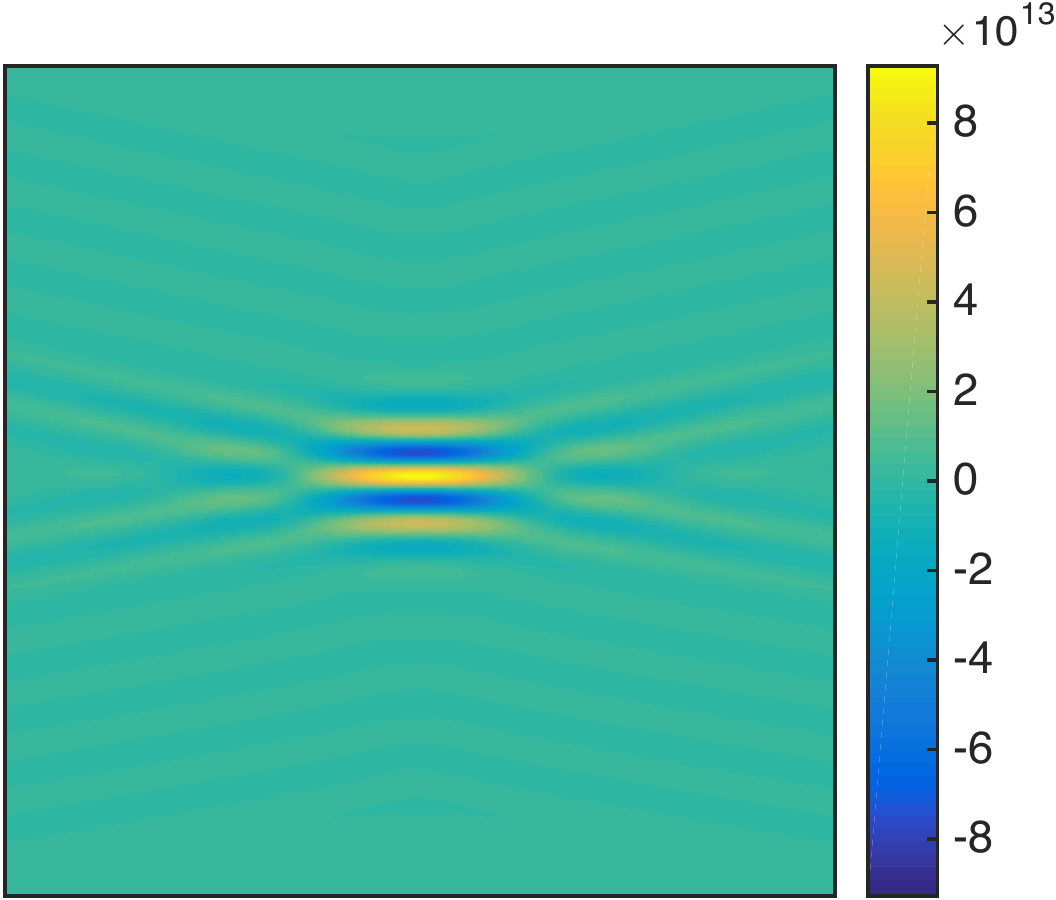}\label{fig:PSF-true-T0}}\hfill{}\subfloat[The PSF with $\theta=0.1$.]{\includegraphics[width=0.24\textwidth]{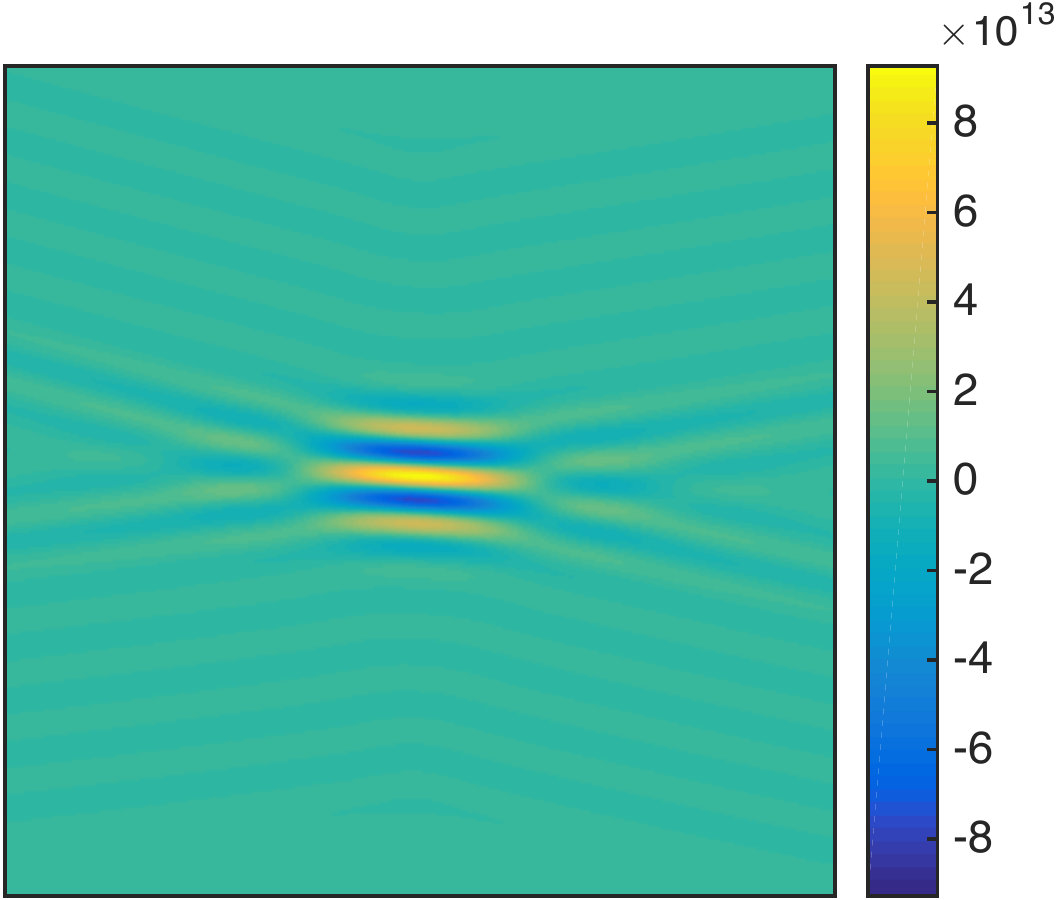}

\label{fig:PSF-true-T01}}\hfill{}\subfloat[The PSF with $\theta=0.2$.]{\includegraphics[width=0.24\textwidth]{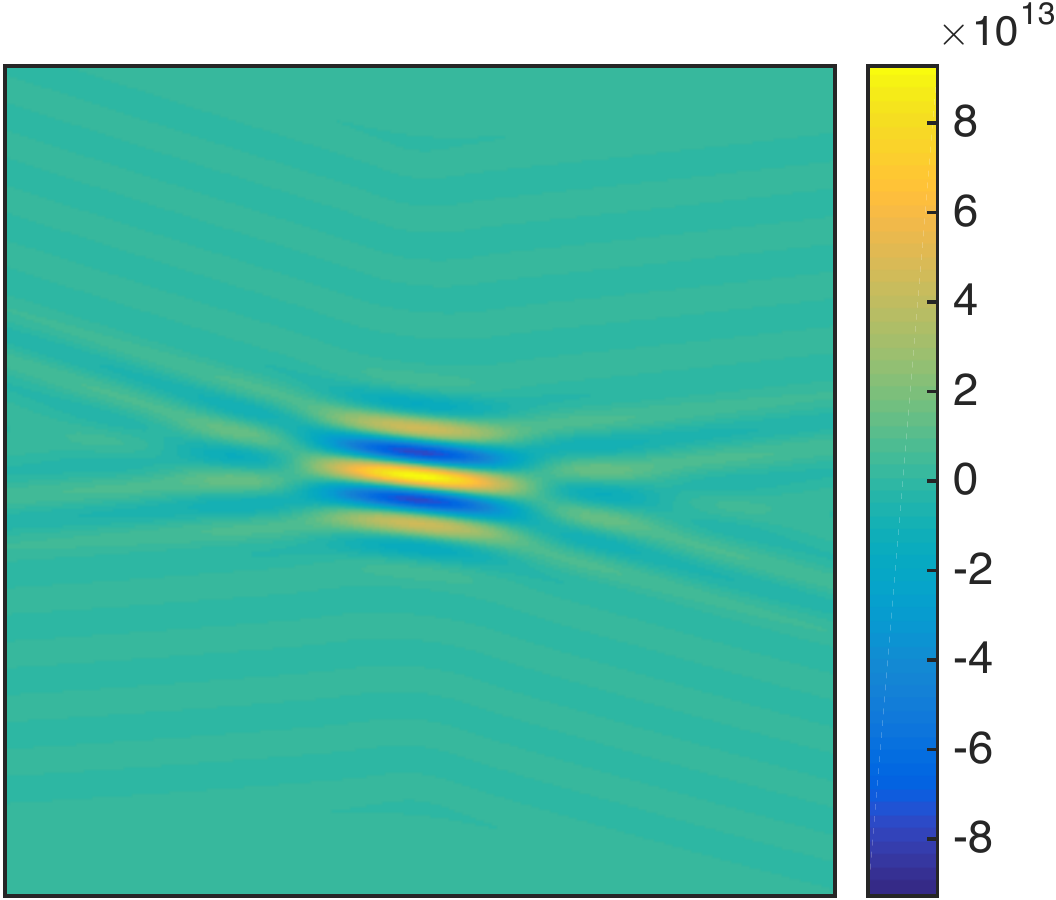}\label{fig:PSF-true-T02}}\hfill{}\subfloat[The PSF with $\theta=0.3$.]{\includegraphics[width=0.24\textwidth]{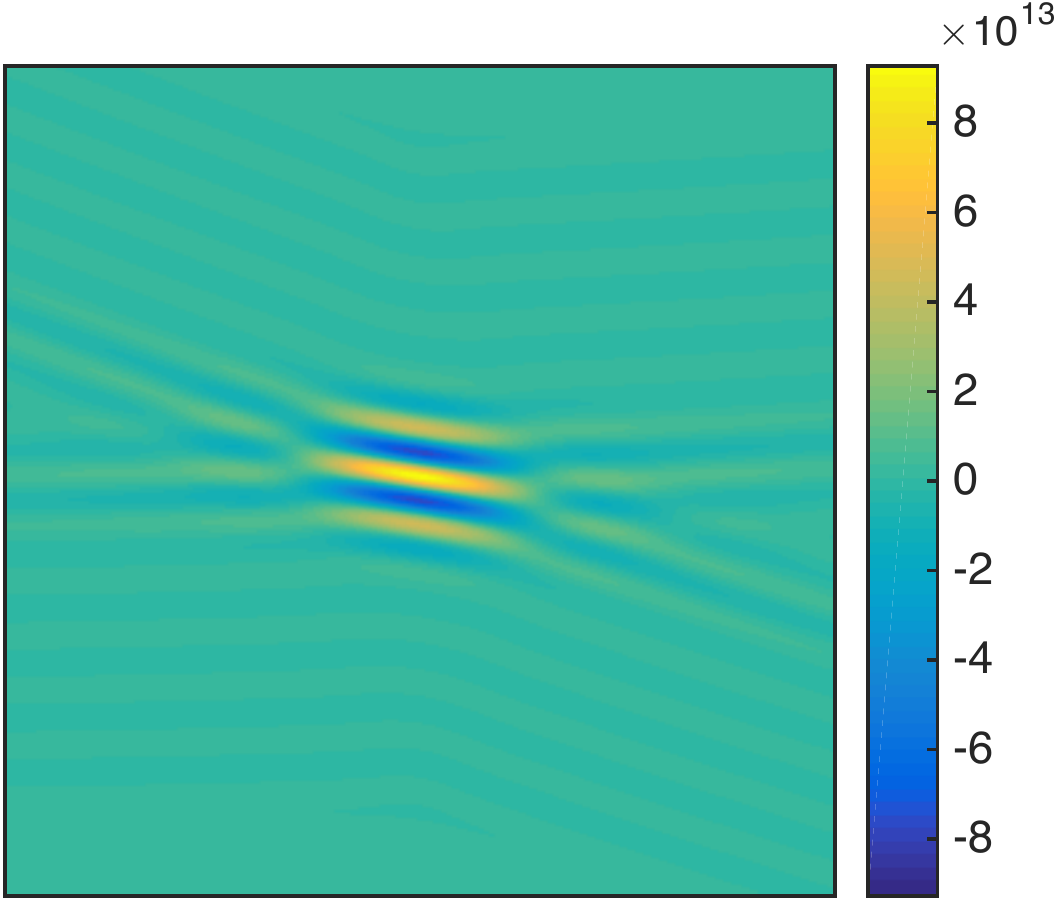}\label{fig:PSF-true-T03}}\caption{The exact PSF with different values of the angle $\theta$  (with parameters as in \eqref{eq:f-chi} and \eqref{eq:quantities}, and $F=0.4$).
The size of the square shown is $\unit[2]{mm}\times\unit[2]{mm}$, and the horizontal and vertical axes are the $x$ and $z$ axes, respectively.}
\label{fig:PSF-T}
\end{figure}

\end{rem}
We have now expressed $g_{\theta}$ as a convolution kernel. In order
to better understand the different roles of the variables $x$ and
$z$, it is instructive to use the actual expression for $f$ given
in (\ref{eq:f-chi}). Since $f'(t)=\nu_{0}e^{2\pi i\nu_{0}t}\tilde{\chi}(\nu_{0}t)$,
with $\tilde{\chi}(t)=2\pi i\chi(t)+\chi'(t)$, we can write
\begin{equation*}
\begin{split}f' & \left(\frac{c_{0}^{-1}}{\sqrt{1+F^{2}}}\left((1+\sqrt{1+F^{2}}\cos\theta)z+(\sqrt{1+F^{2}}\sin\theta\pm F)x\right)\right)\\
 & =\nu_{0}e^{\frac{2\pi i\nu_{0}c_{0}^{-1}}{\sqrt{1+F^{2}}}\left((1+\sqrt{1+F^{2}}\cos\theta)z+(\sqrt{1+F^{2}}\sin\theta\pm F)x\right)}\\
 & \quad\tilde{\chi}\left(\frac{\nu_{0}c_{0}^{-1}}{\sqrt{1+F^{2}}}\left((1+\sqrt{1+F^{2}}\cos\theta)z+(\sqrt{1+F^{2}}\sin\theta\pm F)x\right)\right)\\
 & \approx\nu_{0}e^{2\pi i\nu_{0}c_{0}^{-1}\left(2z+(\theta\pm F)x\right)}\tilde{\chi}\left(2\nu_{0}c_{0}^{-1}z\right),
\end{split}
\end{equation*}
where we have approximated the dependence on $F$ and $\theta$ at
first order around $F=0$ and $\theta=0$ in the complex exponential
(recall that $F$ and $\theta$ are small) and at zero-th order ($F=0$
and $\theta=0$) inside $\tilde{\chi}$: the difference in the orders
is motivated by the fact that the variations of the complex exponentials
have much higher frequencies than those of $\tilde{\chi}$, since several oscillations are contained in the envelope defined by $\chi$, as it can be easily seen in Figure~\ref{fig:pulse} (and similarly for $\chi'$). This
approximation may be justified by arguing as in Appendix~\ref{sec:The-justification-of}.
Inserting this expression into (\ref{eq:g-tilde}) yields
\begin{equation*}
\begin{split}\tilde{g}_{\theta}(\x) & \approx\frac{c_{0}}{4\pi x}\left[\nu_{0}e^{2\pi i\nu_{0}c_{0}^{-1}\left(2z+(\theta-F)x\right)}\tilde{\chi}\left(2\nu_{0}c_{0}^{-1}z\right)\!-\!\nu_{0}e^{2\pi i\nu_{0}c_{0}^{-1}\left(2z+(\theta+F)x\right)}\tilde{\chi}\left(2\nu_{0}c_{0}^{-1}z\right)\right]\\
 & =-\frac{i\nu_{0}c_{0}}{2\pi x}\,\tilde{\chi}\left(2\nu_{0}c_{0}^{-1}z\right)e^{4\pi i\nu_{0}c_{0}^{-1}z}e^{2\pi i\nu_{0}c_{0}^{-1}\theta x}\,\sin(2\pi\nu_{0}c_{0}^{-1}Fx),
\end{split} \end{equation*}
whence for every $\x=(x,z)\in\R^{2}$
\begin{equation}
\tilde{g}_{\theta}(\x)\approx - i\nu_{0}^2F\tilde{\chi}\left(2\nu_{0}c_{0}^{-1}z\right)e^{4\pi i\nu_{0}c_{0}^{-1}z}e^{2\pi i\nu_{0}c_{0}^{-1}\theta x}\,\sinc(2\pi\nu_{0}c_{0}^{-1}Fx),\label{eq:g-sinc}
\end{equation}
where $\sinc(x):=\sin(x)/x$ (see Figure~\ref{fig:PSF-app2}). This final expression allows us
to analyze the PSF $\tilde{g}_{\theta}$, and in particular its different
behaviors with respect to the variables $x$ and $z$. Consider
for simplicity the case $\theta=0$ (with $\tau=1$). In view of the term $\tilde{\chi}\left(2\nu_{0}c_{0}^{-1}z\right)$, the vertical resolution is approximately $0.8\cdot \nu_0^{-1} c_0$; similarly, in view of the term $\sinc(2\pi\nu_{0}c_{0}^{-1}Fx)$, the horizontal resolution is approximately $\frac{1}{2F}\, \nu_0^{-1} c_0$. Even though horizontal and vertical resolutions are comparable, in terms of focusing and frequencies of oscillations the PSF has very different behaviours in the two directions. Indeed, we can observe that the focusing
in the variable $z$ is sharper than that in the variable $x$: the decay of $\tilde{\chi}$
is much stronger than the decay of $\sinc$. Moreover,
in the variable $z$ we have only high oscillations, while in the
variable $x$ the highest oscillations are at least four times slower
($2=4\frac{1}{2}\ge4F$), and very low frequencies are present as
well, due to the presence of the $\sinc$. As it is clear from Figure~\ref{fig:PSF},
this approximation introduces evident distortions of the tails, as
it is expected from the approximation $F=0$ inside $\tilde{\chi}$;
however, the center of the PSF is well approximated. Similar considerations
are valid for the case when $\theta\neq0$: as observed before, this
simply gives a rotation.

\begin{figure}
\captionsetup[subfigure]{width=85pt}\subfloat[The Fourier transform of the exact PSF $g_0$ given in (\ref{eq:def-g}).]{\includegraphics[clip,width=0.23\textwidth]{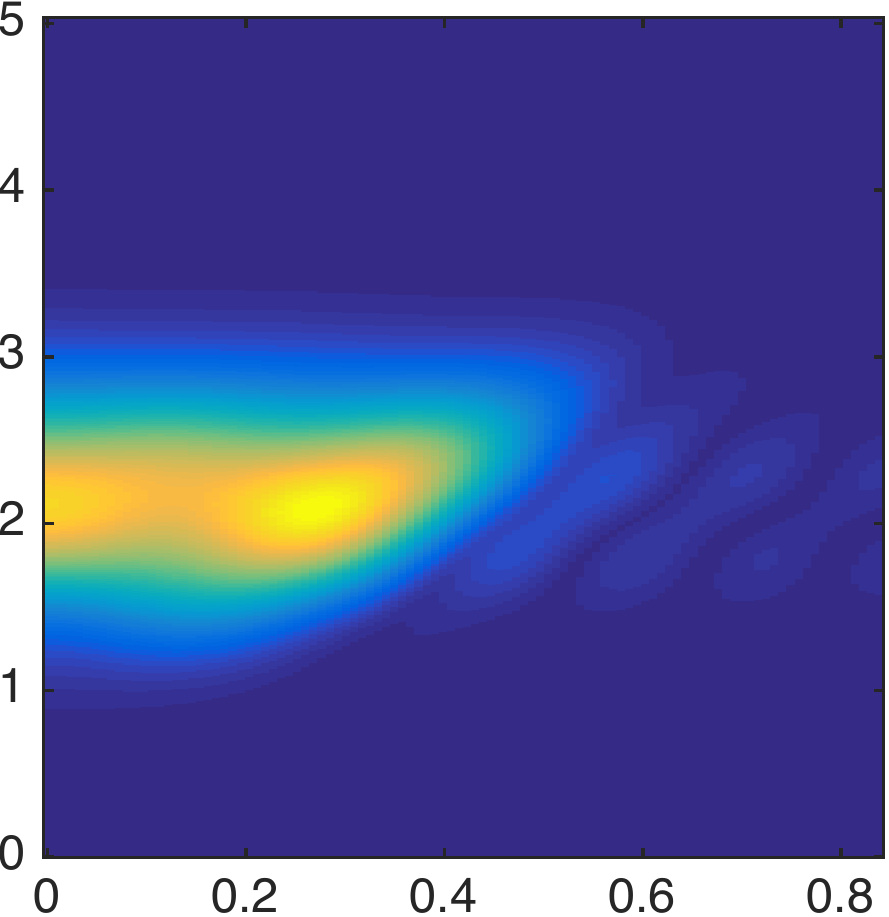}\label{fig:PSF-fourier-true}}\hfill{}\subfloat[The Fourier transform of the approximation of the PSF $\tilde{g}_{0}$
given in (\ref{eq:g-tilde}).]{\includegraphics[width=0.23\textwidth]{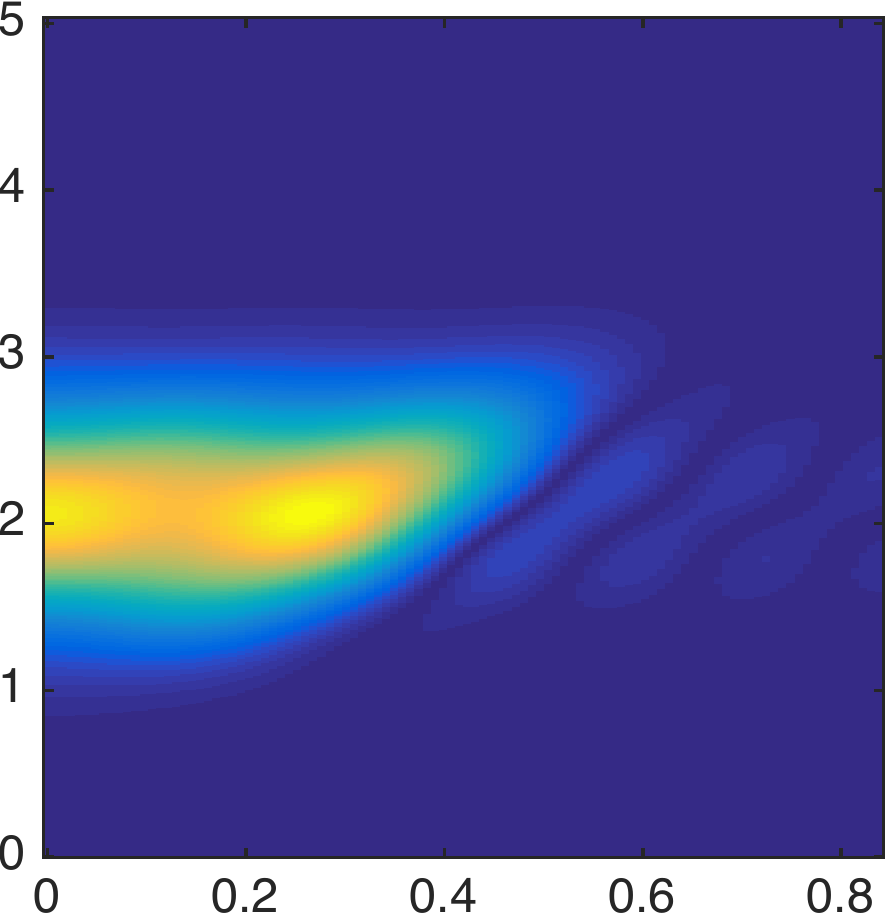}

\label{fig:PSF-fourier-app}}\hfill{}\subfloat[The Fourier transform (\ref{eq:PSF-fourier-sinc}) of the approximation
of the PSF $\tilde{g}_{0}$ given in (\ref{eq:g-sinc}).]{\includegraphics[width=0.23\textwidth]{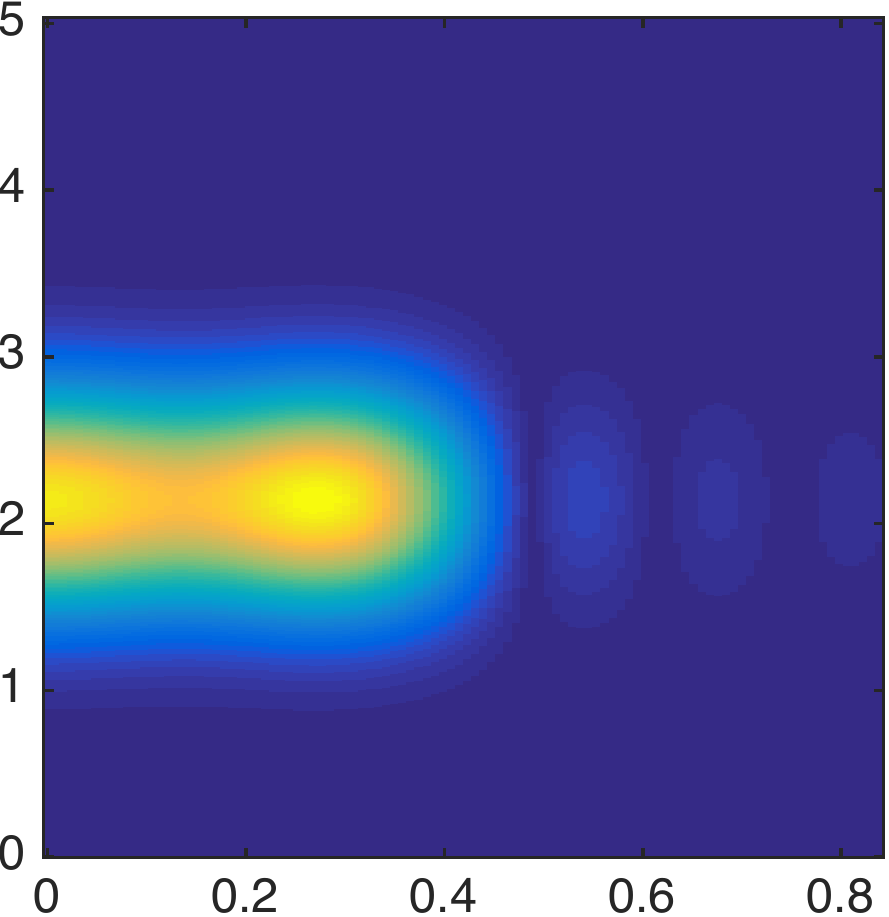}

\label{fig:PSF-fourier-sinc}}\hfill{}\subfloat[The Fourier transform of the PSF $g_{\Theta}^{\text{ac}}$ given in
(\ref{eq:g-angle}), for $\Theta=0.25$.]{\includegraphics[width=0.23\textwidth]{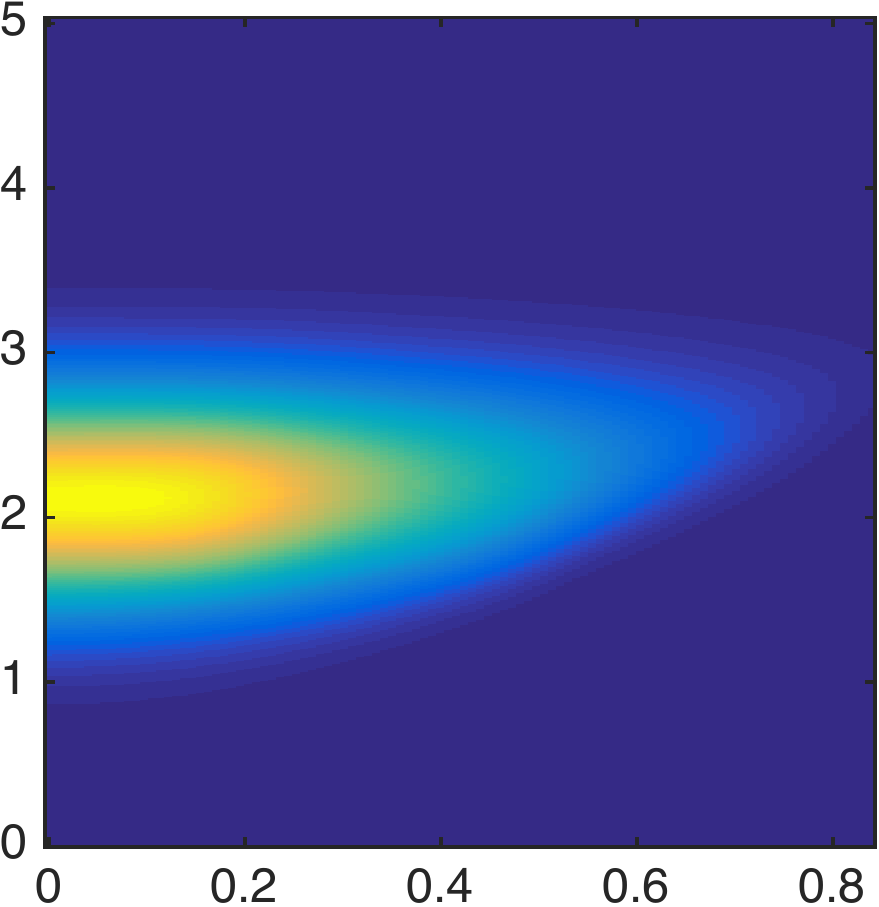}

\label{fig:PSF-fourier-angle}}\caption{The absolute values of the Fourier transforms of the point spread
functions and its approximations (with parameters as in \eqref{eq:f-chi} and \eqref{eq:quantities}, and $F=\theta=0$). The frequency axes are normalized
by $\nu_{0}c_0^{-1}$: the PSF is a low pass filter with cut-off frequency
$F\nu_{0}c_0^{-1}$ with respect to the variable $x$ and a band pass
filter around $2\nu_{0}c_0^{-1}$ with respect to $z$.}
\label{fig:PSF-fourier}
\end{figure}

The same analysis may be carried out by looking at the expression
of the PSF in the frequency domain. For simplicity, consider the case
$\theta=0$: the general case simply involves a translation in the
frequency domain with respect to $x$. Thanks to the separable form
of $\tilde{g}_{\theta}$ given in (\ref{eq:g-sinc}), the Fourier transform may be
directly calculated, and results in the product of the Fourier transform
of $\tilde{\chi}$ and the Fourier transform of the $\sinc$. More
precisely, we readily derive
\begin{equation*}
\begin{split}\F\tilde{g}_{\theta}(\xi_{x},&\xi_{z})  =\int_{\R^{2}}\tilde{g}_{\theta}(x,z)e^{-2\pi i(x\xi_{x}+z\xi_{z})}\,dxdz\\
 & \approx - i\nu_{0}^2F\int_{\R}\sinc(2\pi\nu_{0}c_{0}^{-1}Fx)e^{-2\pi ix\xi_{x}}\,dx\int_{\R}\tilde{\chi}\left(\nu_{0}c_{0}^{-1}z\right)e^{-2\pi i(-2\nu_{0}c_{0}^{-1}+\xi_{z})z}\,dz.
\end{split} \end{equation*}
Thus, since the Fourier transform of the $\sinc$ may be easily computed
and is a suitable scaled version of the rectangle function, we have
\begin{equation*}
\begin{split}\F\tilde{g}_{\theta}(\xi_{x},\xi_{z}) & \approx - i\nu_{0}^2F\frac{1}{2\nu_{0}c_{0}^{-1}F}\mathds{1}_{[-F,F]}\left(c_0\nu_0^{-1}\xi_{x}\right)\int_{\R}\tilde{\chi}\left(\nu_{0}c_{0}^{-1}z\right)e^{-2\pi i(-2\nu_{0}c_{0}^{-1}+\xi_{z})z}\,dz\\
  & =-\frac{ic_{0}\nu_{0}}{2}\mathds{1}_{[-F,F]}\left(c_0\nu_0^{-1}\xi_{x}\right)\frac{1}{\nu_{0}c_{0}^{-1}}\F\tilde{\chi}\left(\frac{-2\nu_{0}c_{0}^{-1}+\xi_{z}}{\nu_{0}c_{0}^{-1}}\right),
\end{split} \end{equation*}
whence
\begin{equation}
\F\tilde{g}_{\theta}(\xi_{x},\xi_{z})\approx - i c_{0}^{2}\,\mathds{1}_{[-F,F]}\left(c_0\nu_0^{-1}\xi_{x}\right)\F\tilde{\chi}\left(-2+\nu_{0}^{-1}c_{0}\xi_{z}\right)/2.\label{eq:PSF-fourier-sinc}
\end{equation}
Therefore, up to a constant, the Fourier transform of the PSF is a
low-pass filter in the variable $x$ with cut-off frequency $F\nu_{0}c_{0}^{-1}$
and a band pass filter in $z$ around $2\nu_{0}c_{0}^{-1}$ (since
$\tilde{\chi}$ is a low-pass filter). This explains, from another
point of view, the different behaviors of $\tilde{g}_{\theta}$ with
respect to $x$ and $z$. This difference is evident from Figure~\ref{fig:PSF-fourier},
where the absolute values of the Fourier transforms of the different
approximations of the PSF are shown.

\subsection{Angle compounding\label{sub:Angle-compounding}}

We saw in the previous subsection that, while very focused in the
direction $z$, the PSF is not very focused in the direction $x$
due to the presence of the $\sinc$ function, see (\ref{eq:g-sinc}).
In order to have a better focusing, it was proposed in \cite{montaldo-tanter-bercoff-benech-finck-2009}
to use multiple measurements corresponding to many angles in an interval
$\theta\in[-\Theta,\Theta]$ for some $0\le \Theta\le0.25$. The reason
why this technique is promising is evident from Figure~\ref{fig:PSF-T}:
adding up several angles together will result in an enhancement of
the center of the PSF, and in a substantial reduction of the artifacts
caused by the tails in the direction $x$. Let us now analyze this
phenomenon analytically.

In a continuous setting, angle compounding corresponds to setting
\begin{equation}
s_{\Theta}^{\text{ac}}(\x)=\frac{1}{2\Theta}\int_{-\Theta}^{\Theta}s_{\theta}(\x)\,d\theta,\qquad\x\in\R_{+}^{2}.\label{eq:s-ac}
\end{equation}
Thus, by linearity, the corresponding PSF is given by
\begin{equation}
g_{\Theta}^{\text{ac}}(\x,\x')=\frac{1}{2\Theta}\int_{-\Theta}^{\Theta}g_{\theta}(\x,\x')\,d\theta,\qquad\x,\x'\in\R_{+}^{2}.\label{eq:g-angle}
\end{equation}

\begin{figure}
\captionsetup[subfigure]{width=110pt}\subfloat[The PSF $g_{\theta}$ with $\theta=0$.]{\includegraphics[clip,width=0.3\textwidth]{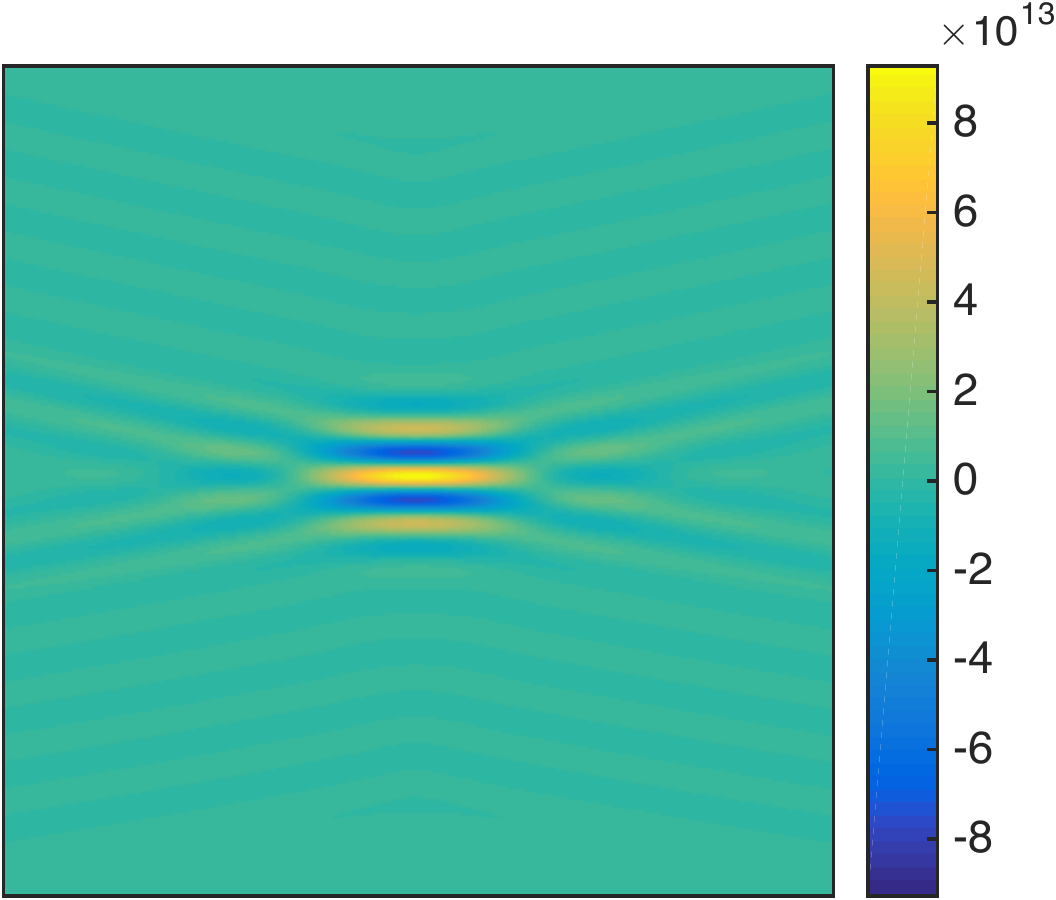}}\hfill{}\subfloat[The PSF $g_{\Theta}^{\text{ac}}$ with $\Theta=0.25$.]{\includegraphics[width=0.3\textwidth]{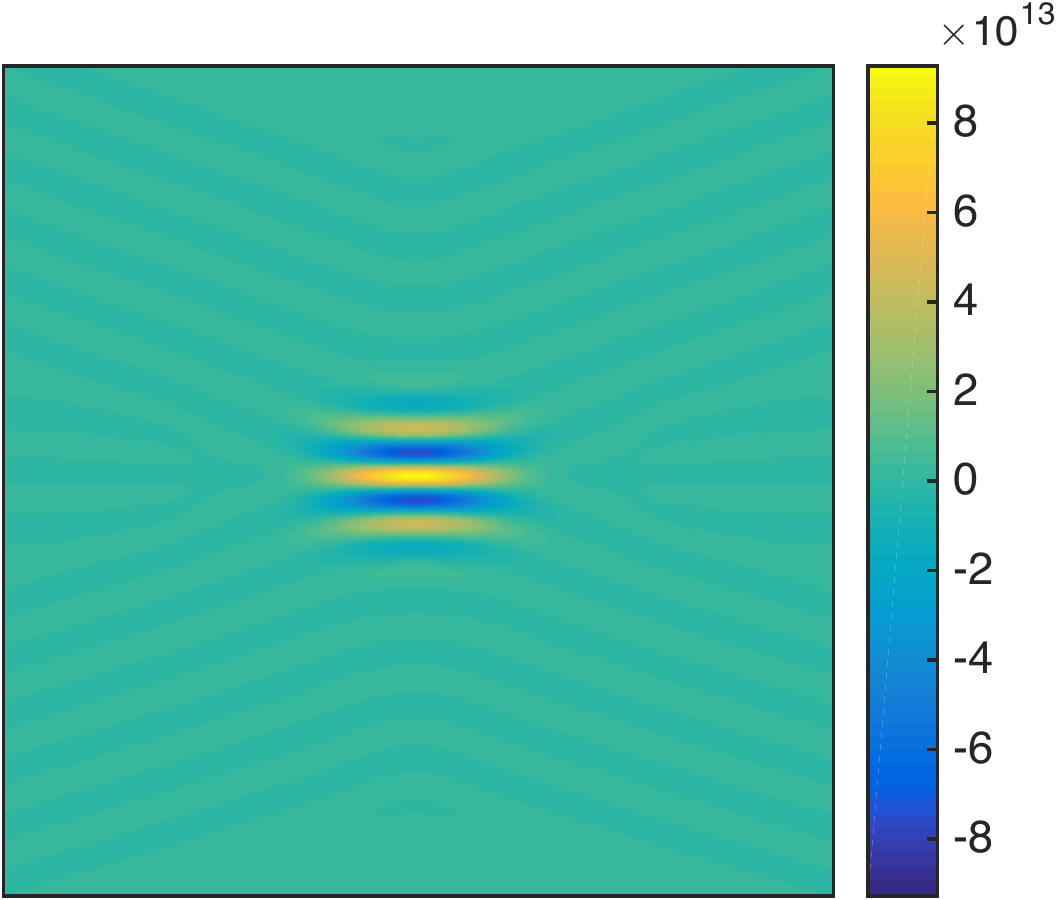}

}\hfill{}\subfloat[The PSF $\tilde{g}_{\Theta}^{\text{ac}}$ with $\Theta=0.25$.]{\includegraphics[width=0.3\textwidth]{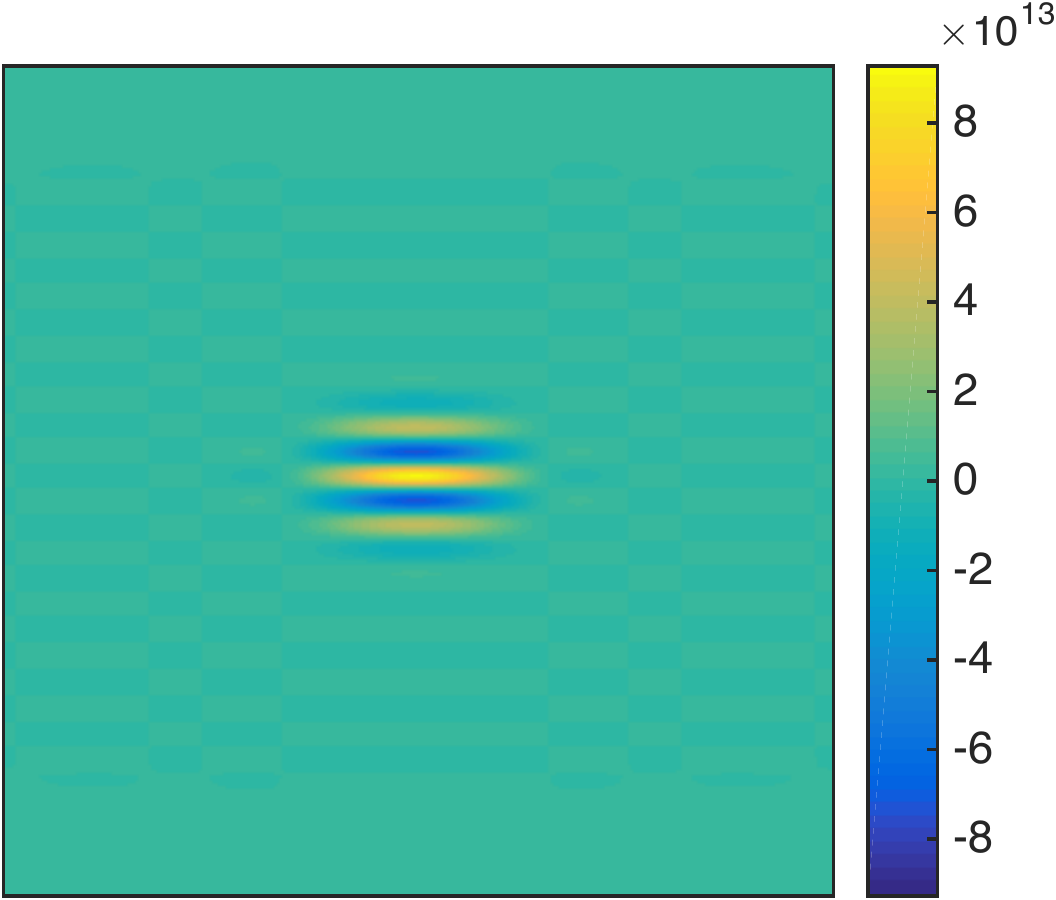}}\caption{A comparison of the PSF related to the single illumination with the
PSF associated to multiple angles (with parameters as in \eqref{eq:f-chi} and \eqref{eq:quantities}, and $F=0.4$). The better focusing in the variable
$x$ for $g_{\Theta}^{\text{ac}}$ is evident, as well as the good
approximation given by $\tilde{g}_{\Theta}^{\text{ac}}$. The size of the square shown is $\unit[2]{mm}\times\unit[2]{mm}$, and the horizontal and vertical axes are the $x$ and $z$ axes, respectively.}
\label{fig:PSF-angle}
\end{figure}

Let us find a simple expression for $g_{\Theta}^{\text{ac}}$.
By using (\ref{eq:g-conv}), we may write $g_{\Theta}^{\text{ac}}(\x,\x')\approx\tilde{g}_{\Theta}^{\text{ac}}(\x-\x')$,
where $\tilde{g}_{\Theta}^{\text{ac}}$ is given by
$
\tilde{g}_{\Theta}^{\text{ac}}(\x)=\frac{1}{2\Theta}\int_{-\Theta}^{\Theta}\tilde{g}_{\theta}(\x)\,d\theta
$, 
so that the image may be expressed as
\begin{equation}
s_{\Theta}^{\text{ac}}(\x)=(\tilde{g}_{\Theta}^{\text{ac}}*n)(\x),\qquad\x\in\R_{+}^{2}.\label{eq:s-ac-approx}
\end{equation}
Thus, in view of the approximation (\ref{eq:g-sinc}), we can write
\begin{equation*}
\begin{split}\tilde{g}_{\Theta}^{\text{ac}}(\x) & =-\frac{i\nu_{0}^2F}{2\Theta}\int_{-\Theta}^{\Theta}\tilde{\chi}\left(2\nu_{0}c_{0}^{-1}z\right)e^{4\pi i\nu_{0}c_{0}^{-1}z}e^{2\pi i\nu_{0}c_{0}^{-1}\theta x}\,\sinc(2\pi\nu_{0}c_{0}^{-1}Fx)\,d\theta\\
 & =-i\nu_{0}^2F\tilde{\chi}\left(2\nu_{0}c_{0}^{-1}z\right)e^{2i\nu_{0}c_{0}^{-1}z}\sinc(2\pi\nu_{0}c_{0}^{-1}Fx)\sinc(2\pi\nu_{0}c_{0}^{-1}\Theta x).
\end{split} \end{equation*}
Therefore, we immediately obtain
\begin{equation}
\tilde{g}_{\Theta}^{\text{ac}}(\x)=\tilde{g}_{0}(\x)\sinc(2\pi\nu_{0}c_{0}^{-1}\Theta x),\qquad\x\in\R^{2}.\label{eq:gac}
\end{equation}
This expression shows that the PSF related to angle compounding is
nothing else than the PSF related to the single angle imaging with
$\theta=0$ multiplied by $\sinc(2\pi\nu_{0}c_{0}^{-1}\Theta x)$.
Thus, for $\Theta=0$ we recover $\tilde{g}_{\theta}$ for $\theta=0$,
as expected. However, for $\Theta>0$, this PSF enjoys faster decay
in the variable $x$. See Figure~\ref{fig:PSF-angle} for an illustration
of $g_{\Theta}^{\text{ac}}$ and $\tilde{g}_{\Theta}^{\text{ac}}$
and a comparison with $g_{\theta}$ and Figure~\ref{fig:PSF-fourier-angle}
for an illustration of the Fourier transform of $g_{\Theta}^{\text{ac}}$.

To sum up the main features of the static problem, we have shown that
the recovered image may be written as $s_{\Theta}^{\text{ac}}=\tilde{g}_{\Theta}^{\text{ac}}*n$,
where $\tilde{g}_{\Theta}^{\text{ac}}$ is the PSF of the imaging
system with measurements taken at multiple angles. The ultrafast imaging
technique is based on obtaining many of these images over time, as
we discuss in the next section.

\section{The Dynamic Forward Problem}  \label{sect3}

\subsection{The quasi-static approximation and the construction of the data}

The dynamic imaging setup consists in the repetition of the static
imaging method over time to acquire a collection of images of a medium
in motion. We consider a quasi-static model: the whole process of
obtaining one image, using the image compounding technique discussed
in Subsection~\ref{sub:Angle-compounding}, is fast enough to consider
the medium static, but collecting several images over time gives us
a movie of the movement over time. In other words, there
are two time scales: the fast one related to the propagation of the
wave is considered instantaneous with respect to the slow one, related
to the sequence of the images.

In view of this quasi-static approximation, from now on we neglect
the time of the propagation of a single wave to obtain static imaging.
The time $t$ considered here is related to the slow time scale. In
other words, by (\ref{eq:s-ac-approx}) at fixed time $t$ we obtain
a static image $s(\mathbf{x},t)$ of the medium $n=n(\mathbf{x},t)$,
namely
\begin{equation}\label{eq:dynamic-s}
s(\mathbf{x},t)=\left(\tilde{g}_{\Theta}^{\text{ac}}*n(\,\cdot\,,t)\right)(\mathbf{x}).
\end{equation}
Repeating the process for $t\in[0,T]$ we obtain the movie $s(\mathbf{x},t)$,
which represents the main data we now need to process. As mentioned
in the introduction, our aim is locating the blood vessels within
the imaged area, by using the fact that $s(\x,t)$ will be strongly
influenced by movements in $n$.

\subsection{The doppler effect\label{sub:doppler}}

Measuring the medium speed is an available criterion to separate different
sources; thus, we want to see the influence on the image of a
single particle in movement, as by linearity the obtained conclusions
 naturally extend to a group of particles. For a single particle,
we are interested in observing the generated Doppler effect in the
reconstructed image,  namely peaks in the Fourier transform away from zero.

Intuitively, Figure \ref{fig:PSF-fourier-angle} shows that there
is a clear difference in the movements depending on their orientation. We
want to explore this difference in a more precise way. Let us consider
$n(x,z,t)=\delta_{(0,vt)}(x,z)$, i.e.\ a single particle moving in the
$z$ direction with velocity $v$. The resulting image, as a function
of time, is obtained via equations  (\ref{eq:gac}) and \eqref{eq:dynamic-s}
\begin{equation*}
\begin{split}
s(x,z,t)  &\approx  \int_{\R^{2}}\tilde{g}_{\Theta}^{\text{ac}}(x-x',z-z')\delta_{(0,vt)}(x',z')dx'dz'\\
& =  \tilde{g}_{\Theta}^{\text{ac}}(x,z-vt)\\
 & =  \tilde{g}_{0}(x,z-vt)\sinc(2\pi\nu_{0}c_{0}^{-1}\Theta x).
\end{split} \end{equation*}
Therefore, arguing as in (\ref{eq:PSF-fourier-sinc}), we obtain that the Fourier transform with respect to the time variable $t$ of the image is given by
\begin{equation*}
\begin{split}
\F_{t}(s)(x,z,\xi) & \approx  \int_{\R}\tilde{g}_{0}(x,z-vt)e^{-2\pi i\xi t}dt\ \sinc(2\pi\nu_{0}c_{0}^{-1}\Theta x)\\
 & =  \frac{1}{v}e^{-2\pi i\frac{\xi z}{v}}\F_{2}(\tilde{g}_{0})(x,-\frac{\xi}{v})\sinc(2\pi\nu_{0}c_{0}^{-1}\Theta x),
\end{split} \end{equation*}
where $\F_{2}$ is the Fourier transform with respect to the  variable $z$. Adopting approximation (\ref{eq:g-sinc}), we
obtain
\[
\F_{t}(s)(x,z,\xi)  \!\approx\! -  \frac{1}{v}i\nu_{0}^2Fe^{-2\pi i\frac{\xi z}{v}}\sinc(2\pi\nu_{0}c_{0}^{-1}\Theta x)\sinc(2\pi\nu_{0}c_{0}^{-1}Fx)\F(\tilde{\chi})\Bigl(\frac{-\xi}{2\nu_{0}c_{0}^{-1}v}-1\Bigr).
\]
 Given the shape of $\tilde{\chi}$, its Fourier transform has a maximum
around 0, thus we can see a peak of $|\F_{t}(s)(x,z,\xi)|$ when $\xi$
is around $-2\nu_{0}c_{0}^{-1}v$, and so we have the Doppler effect.

In the case when the particle is moving parallel to the detector
array, namely $n(x,z,t)=\delta_{(vt,0)}(x,z)$, following an analogous
procedure as before, we obtain
\[
s(x,z,t)  \approx  \tilde{g}_{0}(x-vt,z)\sinc(2\pi\nu_{0}c_{0}^{-1}\Theta(x-vt)),
\]
and applying the Fourier transform in time yields
\[
\F_{t}(s)(x,z,\xi)  \approx  \frac{1}{v}e^{-2\pi i\frac{\xi x}{v}}\F(\tilde{g}_{0}(\cdot,z)\sinc(2\pi\nu_{0}c_{0}^{-1}\Theta\cdot))(-\frac{\xi}{v}).
\]
Using approximation (\ref{eq:g-sinc}), the convolution formula for
the Fourier transform and the known transform of the $\sinc$ function,
gives
\begin{equation*}
\F_{t}(s)(x,z,\xi)  \approx - i\frac{e^{-2\pi i\frac{\xi x}{v}}}{4\Theta v}\nu_0 c_{0}\tilde{\chi}(2\nu_{0}c_{0}^{-1}z)e^{4\pi i\nu_{0}c_{0}^{-1}z}(\mathds{1}_{[-F,F]}*\mathds{1}_{[-\Theta,\Theta]})\left(-\frac{\xi}{v\nu_{0}c_{0}^{-1}}\right).
\end{equation*}
The convolution of these characteristic functions evaluated at $\eta$ is equal
to the length of interval $[-F+\eta,F+\eta]\cap[-\Theta,\Theta]$,
because
\[
(\mathds{1}_{[-F,F]}*\mathds{1}_{[-\Theta,\Theta]})(\eta)  \!=\!  \int_{\R}\mathds{1}_{[-F,F]}(\eta-s)\mathds{1}_{[-\Theta,\Theta]}(s)ds\! =\! \int_{\R}\mathds{1}_{[-F+\eta,F+\eta]}(s)\mathds{1}_{[-\Theta,\Theta]}(s)ds.
\]
Since both intervals are centered at $0$, this value is maximized
for $\eta$ (and thus $\xi$) around $0$, like in the static case, and so the observed
Doppler effect is very small.

These differences are fundamental to
understand  the capabilities of the method for blood flow imaging. This phenomenon will be  experimentally verified in Section~\ref{sect5}.

\subsection{Multiple scatterer random model\label{sub:probmodel}}

We have seen the effect on the image $s(x,z,t)$ of a single moving particle. We now consider the more realistic case of a medium (either blood vessels or tissue) with a
large number of particles in motion. This will allow to study the
statistical properties of the resulting measurements.

We consider a rectangular domain  $\Omega=(-L_{x}/2,L_{x}/2)\times(0,L_{z})$, which consists in $N$
point particles. Let us denote the location of particle $k$ at time
$t$ by $a_{k}(t)$. In the most general case, each particle is subject
to a dynamics
\begin{equation}
a_{k}(t)=\varphi_{k}\left(u_{k},t\right),\qquad a_{k}(0)=u_{k},\label{eq:a_k}
\end{equation}
where $\left(u_{k}\right)_{k=1,...,N}$ are independent uniform random
variables on $\Omega$ and $\left(\varphi_{k}\right)_{k=1,...,N}$
are independent and identically distributed stochastic flows: for instance, they can
 be the flows of a stochastic differential equation or 
the deterministic flows of a partial differential equation.    Thus, the $a_{k}$s are independent and identically
distributed stochastic processes. In view of these considerations, we consider the medium given by
\begin{equation}\label{eq:n_model-0}
n\left(\mathbf{x},t\right)=\frac{C}{\sqrt{N}}\sum_{k=1}^{N}\delta_{a_{k}\left(t\right)}\left(\mathbf{x}\right),
\end{equation}
where $C>0$ denotes the scattering intensity and  $\frac{1}{\sqrt{N}}$ is the natural normalization factor in view of 
the central limit theorem.

To avoid minor issues from boundary effects, which are of no interest
to us in the analysis of this problem, we assume the periodicity of
the medium. In other words, we consider the periodization
\begin{equation}
n_{p}(\x,t)=\sum_{\mathbf{l}\in\Z^{2}}n(\x+\mathbf{l}\cdot\mathbf{L},t),\label{eq:n_model}
\end{equation}
where $\mathbf{L}=(L_{x},L_{z})$.
Let $g\left(\x\right):=\sum_{\mathbf{l}\in\Z^{2}}\tilde{g}_{\Theta}^{\text{ac}}\left(\x+\mathbf{l}\cdot\mathbf{L}\right)$
be the periodic PSF, which is more convenient than $\tilde{g}_{\Theta}^{\text{ac}}$
(given by (\ref{eq:gac})) for a $\Omega$-periodic medium. The dynamic image
$s$ is then given by
\[
s(\mathbf{x},t)=(\tilde{g}_{\Theta}^{\text{ac}}*n_{p}\left(\cdot,t\right))\left(\mathbf{x}\right)=(g*n(\,\cdot\,,t))(\x)=\frac{C}{\sqrt{N}}\sum_{k=1}^{N}g\left(\mathbf{x}-a_{k}\left(t\right)\right).
\]
Let us also
assume for the sake of simplicity that,
at every time $t$,  $a_{k}\left(t\right)$ modulo  $\Omega$ is a uniform random variable
on $\Omega$, namely
\begin{equation}\label{eq:uniform}
\E \sum_{\mathbf{l}\in\Z^2} w(a_k(t)+\mathbf{l}\cdot \mathbf{L}) =  |\Omega|^{-1}\int_{\R^2} w(\mathbf{y})\,d\mathbf{y},\qquad w\in L^1(\R^2).
\end{equation}
As a simple but quite general example, it is worth noting that in the case when $a_k(t) = u_k + F(t)$, where $F(t)$ is any random process independent of $u_k$, the above equality is satisfied, since
\[
\E \sum_{\mathbf{l}\in\Z^2} w(u_k + F(t) +\mathbf{l}\cdot \mathbf{L}) = |\Omega|^{-1}\E \sum_{\mathbf{l}\in\Z^2} \int_\Omega w(\mathbf{y} + F(t) +\mathbf{l}\cdot \mathbf{L}) d\mathbf{y}  = |\Omega|^{-1} \int_{\R^2} w(\mathbf{y})\,d\mathbf{y},
\]
where the expectation in the first term is taken with respect to $u_k$ and $F(t)$, while in the second term only with respect to $F(t)$.

We now wish to compute the expectation of the random variables present in the expression for $s(\x,t)$. By \eqref{eq:g-sinc} and (\ref{eq:gac}), since $\tilde{g}_{\Theta}^{\text{ac}}$ is a derivative of a Schwartz function in the variable $z$, we have $\int_{\R^2}\tilde{g}_{\Theta}^{\text{ac}}(\mathbf{y})d\mathbf{y}=0$. Thus, by \eqref{eq:uniform}
the expected value may be easily computed as
\begin{equation}
\E\left(g\left(\mathbf{x}-a_{k}\left(t\right)\right)\right)= \E \sum_{\mathbf{l}\in\Z^2} \tilde{g}_{\Theta}^{\text{ac}}(\x-a_k(t)+\mathbf{l}\cdot \mathbf{L}) = |\Omega|^{-1} \int_{\R^2} \tilde{g}_{\Theta}^{\text{ac}} (\mathbf{y})d\mathbf{y} = 0.\label{eq:expectation}
\end{equation}

Let $\left(\mathbf{x}_{i}\right)_{i=1,...,m_{\x}}$ and $\left(t_{j}\right)_{j=1,...,m_{t}}$
be the sampling locations and times respectively. The data may be
collected in the Casorati matrix $S_{N}\in\mathbb{C}^{m_{\x}\times m_{t}}$
defined by
\[
S_{N}(i,j)=s(\mathbf{x}_{i},t_{j}).
\]
By (\ref{eq:expectation}), according to the multivariate central
limit theorem, the matrix $S_{N}$ converges in distribution to a
Gaussian complex matrix $S\in\mathbb{C}^{m_{\x}\times m_{t}}$ , the
distribution of which is entirely determined by the following correlations,
for $i,i'=1,\dots,m_{\x}$ and $j,j'=1,\dots,m_{t}$
\begin{align}
 & \E(S(i,j))=0, \nonumber \\
 & \Cov\left(S(i,j),S(i',j')\right)=C^{2}\E\left(g\left(\mathbf{x}_{i}-a_{1}\left(t_{j}\right)\right)g\left(\mathbf{x}_{i'}-a_{1}\left(t_{j'}\right)\right)\right),\label{eq:ss}\\
 & \Cov\left(S(i,j),\overline{S(i',j')}\right)=C^{2}\E\left(g\left(\mathbf{x}_{i}-a_{1}\left(t_{j}\right)\right)\overline{g\left(\mathbf{x}_{i'}-a_{1}\left(t_{j'}\right)\right)}\right).\label{eq:ssbar}
\end{align}
More precisely, let $w\in\C^{m_{\x}m_{t}}$ be a column vector containing
all the entries of $S$. Let $v\in\C^{2m_{\x}m_{t}}$ and $V\in\C^{2m_{\x}m_{t}\times2m_{\x}m_{t}}$
be defined by
\[
v=\left(w_{1},\overline{w_{1}},w_{2},\overline{w_{2}},...,w_{m_{\x}m_{t}},\overline{w_{m_{\x}m_{t}}}\right)^{T}\quad\text{and}\quad V=\E\left(v\overline{v}^{T}\right).
\]
The covariance matrix $V$ can be easily computed from~(\ref{eq:ss})
and (\ref{eq:ssbar}). Then the probability density function $f$
of $v$ can be expressed as \cite{bos_multivariate_1995}:
\[
f\left(v\right)=\frac{1}{\pi^{m_{\x}m_{t}}\det\left(V\right)^{\frac{1}{2}}}\exp\left(-\frac{1}{2}v^{*}V^{-1}v\right).
\]
Moreover, it is possible to generate samples from this distribution:
if $X$ is a complex unit variance independent normal random vector, and if
$\sqrt{V}$ is a square root of $V$, then $\sqrt{V}X$ is distributed
like $v$. This allows for simulations of sample image sequences for 
a large number of particles with a complexity independent 
of the number of particles.

The analysis carried out here will allow us to study the distribution
of the singular value of the matrix $S$, depending on the properties
of the flows $\varphi_{k}$. This will be the key ingredient to justify
the correct separation of blood and clutter signals by means of the
singular value decomposition of the measurements.

\section{The Dynamic Inverse Problem: Source Separation}  \label{sect4}

\subsection{Formulation of the dynamic inverse problem}

As explained in the introduction, the aim of the dynamic inverse problem
is blood flow imaging. In other words, we are interested in locating
blood vessels, possibly very small, within the medium. The main issue
is that the signal $s(\x,t)$ is highly corrupted by clutter signal,
namely the signal scattered from tissues. In the linearized regime we consider,
we may write the refractive index $n$ as the sum of a clutter component
$n_{c}$ and a blood component $n_{b}$, namely $n=n_{c}+n_{b}$.
Blood is located only in small vessels in the medium, whereas clutter
signal comes from everywhere: by~(\ref{eq:n-tilde}), since
blood vessels are smaller than the focusing height, even pixels
located in blood vessels contain reflections coming from the tissue.
Let us denote the location of blood vessels by $\Omega_{b}\subset\Omega$.
The inverse problem is the following: can we
recover $\Omega_{b}$ from the data $s(\x,t)=s_{c}(\mathbf{x},t)+s_{b}(\mathbf{x},t)$? Here, $s_c$ and $s_b$ are given by \eqref{eq:dynamic-s}, with $n$ replaced by $n_c$ and $n_b$, respectively.
In this section, we provide a quantitative analysis of the method
described in \cite{demene2015spatiotemporal} based on the singular
value decomposition (SVD) of $s$.

\subsection{The SVD algorithm}

We now review the SVD algorithm presented in \cite{demene2015spatiotemporal}.
The Casorati matrix $S\in\C^{m_\x\times m_t}$ is defined as in previous section by
\[
S(i,j)=s\left(\x_{i},t_{j}\right),\qquad i\in\left\{ 1,...,m_\x\right\},\;j\in\left\{ 1,...,m_t\right\} .
\]
Without loss of generality, we further assume that $m_t\le m_\x$. We remind the reader that the SVD of $S$ is given by
\[
S=\sum_{k=1}^{m_t}\sigma_{k} u_{k} \overline{v_{k}}^T,
\]
where $\left(u_{1},...,u_{m_\x}\right)$ and $\left(v_{1},...,v_{m_t}\right)$
are orthonormal bases of $\C^{m_x}$ and $\C^{m_t}$, and $\sigma_{1}\geq\sigma_{2}\geq...\geq\sigma_{m_t}\ge 0$.
For any $K\geq1$, $S_{K}=\sum_{k=1}^{K}\sigma_{k} u_{k} \overline{v_{k}}^T$
is the best rank $K$ approximation of $S$ in the Frobenius norm.
The SVD is a well-known tool for denoising sequences of images,
see for example \cite{iizuka1982data}. The idea is that since singular
values for the clutter signal are quickly decaying after a certain
threshold, the best rank $K$ approximation
of $S$ will contain most of the signal coming from the clutter, provided that $K$ is large enough. This
could be used to recover clutter data, by applying a ``denoising''
algorithm, and keeping only $S_K$. But it can also be used to recover the blood location, by considering the ``power doppler''
\[
\hat{S}_{b,K}\left(i\right):=\sum_{k=K+1}^{m_t}\sigma_{k}^{2}|u_{k}|^{2}(i)=\sum_{j=1}^{m_t}\left|\left(S-S_{K}\right)(i,j)\right|^{2},\qquad i\in\left\{ 1,...,m_\x\right\}.
\]
As we will show in the following subsection, clutter signal can be
well approximated by a low-rank matrix. Therefore, $S_{K}$
will contain most of the clutter signal for $K$  large enough. In
this case, even if the intensity of total blood reflection is small, $S-S_{K}$
will contain more signal coming from the blood than from the clutter and
therefore high values of $\hat{S}_{b,K}\left(i\right)$ should
be located in blood vessels. 

Before presenting the justification of this method, let us briefly provide a heuristic motivation by considering the SVD of the continuous data given by
\[
s(\x,t)=\sum_{k=1}^\infty \sigma_k u_k(\x)\overline{v_k}(t).
\]
In other words, the dynamic data $s$ is expressed as a sum of spatial components $u_k$ moving with time profiles $\overline{v_k}$, with weights $\sigma_k$. Therefore, since the tissue movement  has higher spatial coherence than the blood flow, we expect the first factors to contain the clutter signal, and the remainder to provide information about the blood location via the quantity $\hat{S}_{b,K}$.

\subsection{Justification of the SVD in 1D}

We will assume that the particles of the blood and of the clutter have independent
dynamics described by (\ref{eq:a_k})-(\ref{eq:n_model}). We add the subscripts
$b$ and $c$ to indicate the dynamics of blood and clutter, respectively.

In this subsection, using the limit Gaussian model presented in $\S$\ref{sub:probmodel}, we present the statistics of the singular values in a
simple 1D model. These are useful to understand the behavior of SVD filtering.
The results of $\S$\ref{sub:probmodel} allow to simulate large number of sample signals $s$, given that we can compute
the covariance matrices~(\ref{eq:ss}) and (\ref{eq:ssbar}). Since
these matrices are very large, we restrict ourselves to the $1D$ case,
so that all sampling locations $\mathbf{x}_{i}$ are located at $x=0$,
and are thus characterized by their depth $z_{i}$. We will therefore
drop all references to $x$ in the following. We also consider very
simplified dynamics, which can be thought of as local descriptions of
the global dynamics at work in the medium. Let $a_{b}=a_{1,b}$ and
$a_{c}=a_{1,c}$ be the random variables for the dynamics of blood
and clutter particles, respectively, as introduced in (\ref{eq:a_k}).
The dynamics is modelled by a Brownian motion with drift, namely
\[
a_{\alpha}\left(t\right)=u_{\alpha}+v_{\alpha}t+\sigma_{\alpha}B_{t},\qquad\alpha\in\left\{ b,c\right\} .
\]
Here, $u_{\alpha}$ represents the position of the particle at time
$t=0$, and is uniformly distributed in $(0,L_{\alpha})$, where $L_{b}\ll L_{c}$.
The deterministic quantity $v_{\alpha}$ is the mean velocity of the
particles. In order to take into account the random fluctuations of
the particles in movement, we added a diffusion term $\sigma_{\alpha}B_{t}$,
where $B_{t}$ is a Brownian motion and $\sigma_{\alpha}^{2}$ is
a diffusion coefficient quantifying the variance of the fluctuations
of the particle position relative to the mean trajectory. We also make the 
simplifying assumption that the diffusion terms are independent over different particles.
More precisely,
we have the following conditional expectation and variance:
\[
\E\left(\left.a_{\alpha}\left(t\right)\right|u_{\alpha}\right)=u_{\alpha}+v_{\alpha}t,\qquad\Var\left(\left.a_{\alpha}\left(t\right)\right|u_{\alpha}\right)=t\sigma_{\alpha}^{2}.
\]

The difference between clutter and blood dynamics is in the diffusion
coefficient: in the case of clutter, since it is an elastic displacement,
$\sigma_{c}^{2}\approx0$. For simplicity, from now on we set $\sigma_{c}=0$.
In the case of blood, which is modelled as a suspension of cells in
a fluid, we have $\sigma_{b}^{2}=\sigma^{2}>0$. This coefficient is expressed
in $\unit{m^{2}s^{-1}}$, and models the random diffusion in a fluid transporting
red blood cells due to turbulence in the fluid dynamics and collisions
between cells. In practice, $\sigma^{2}$ is much larger
than the diffusion coefficient of microscopic particles in a static
fluid, and depends on the velocity $v_{b}$ \cite{caro2012mechanics}.
As for the mean velocities, in the most extreme cases, $v_{b}$ and
$v_{c}$ can be of the same order, even though most of the time $v_{b}>v_{c}$.

Let $S_{b}$ and $S_{c}$ denote the data matrix constructed in $\S$\ref{sub:probmodel},
related to blood and clutter signal, respectively. We now compute
the covariance matrix $V$ of $S_{\alpha}$:
\begin{equation*}
\begin{split}\Cov(S_{\alpha}(i,j&),S_{\alpha}(i',j'))  =C_{\alpha}^{2}\E\left(g\left(z_{i}-a_{\alpha}\left(t_{j}\right)\right)g\left(z_{i'}-a_{\alpha}\left(t_{j'}\right)\right)\right)\\
 & =\frac{C_{\alpha}^{2}}{L}\E\!\int_{0}^{L}g\left(z_{i}-y-v_{\alpha}t_{j}-\sigma_{\alpha}v_{\alpha}B_{t_{j}}\right)g\bigl(z_{i'}-y-v_{\alpha}t_{j'}-\sigma_{\alpha}v_{\alpha}B_{t_{j'}}\bigr)dy\\
 & =C_{\alpha}^{2}\E\,C_{gg}\Bigl(z_{i}-z_{i'}+v_{\alpha}\bigl(t_{j'}-t_{j}+\sigma_{\alpha}(B_{t_{j'}}-B_{t_{j}})\bigr)\Bigr),
\end{split} \end{equation*}
where $C_{gg}\left(z\right)=\frac{1}{L}\int_{0}^{L}g(y)g(z+y)dy$
and $C_{b}$ and $C_{c}$ denote the intensity of the blood and clutter
signals, respectively. The expectation operator is taken over all possible positions $u_\alpha$ and all possible drifts $B_{t_{j}}$ and $B_{t_{j'}}$  in the first line, and only over all drifts in the second and third lines. By standard properties of the Brownian motion,
$B_{t_{j'}}-B_{t_{j}}$ is Gaussian distributed, of expected value
$0$ and variance $\left|t_{j}-t_{j'}\right|$ and so it has the same
distribution as $B_{t_{j'}-t_{j}}$. Thus, in the case of the blood,
we can write
\[
\Cov\left(S_{b}(i,j),S_{b}(i',j')\right)=C_{b}^{2}\E\,C_{gg}\bigl(z_{i}-z_{i'}+v_{b}(t_{j'}-t_{j}+\sigma_{b}B_{t_{j'}-t_{j}})\bigr).
\]
Likewise,
\[
\Cov\bigl(S_{b}(i,j),\overline{S_{b}(i',j')}\bigr)=C_{b}^{2}\E\,C_{g\bar{g}}\bigl(z_{i}-z_{i'}+v_{b}(t_{j'}-t_{j}+\sigma_{b}B_{t_{j'}-t_{j}})\bigr),
\]
where $C_{g\bar{g}}\left(z\right)=\frac{1}{L}\int_{0}^{L}g(y)\bar{g}(z+y)dy$.
The tissue model is then given by $\sigma_{c}=0$, and is therefore
deterministic given the initial position. Thus
\begin{align*}
&\Cov\left(S_{c}(i,j),S_{c}(i',j')\right)=C_{c}^{2}C_{gg}\left(z_{i}-z_{i'}+v_{c}\left(t_{j'}-t_{j}\right)\right),\\
&\Cov\bigl(S_{c}(i,j),\overline{S_{c}(i',j')}\bigr)=C_{c}^{2}C_{g\bar{g}}\left(z_{i}-z_{i'}+v_{c}\left(t_{j'}-t_{j}\right)\right).
\end{align*}
\begin{figure}
\centering
\subfloat[The real part of  $C_{gg}$\label{fig:cgg}.]{\includegraphics[width=0.45\textwidth]{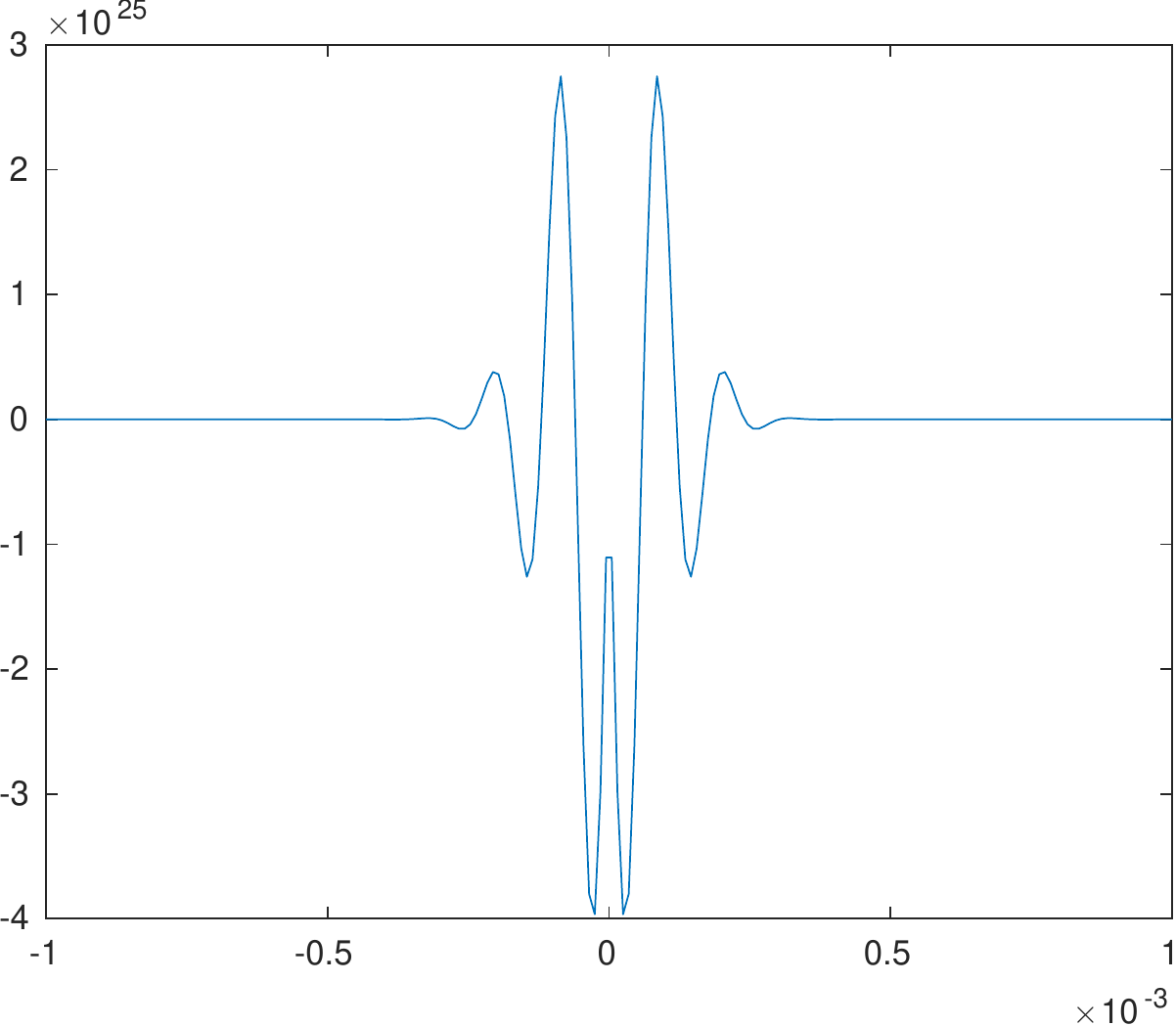}
}\qquad\subfloat[The real part of  $C_{g\overline{g}}$.\label{fig:cgg2}]{\includegraphics[width=0.45\textwidth]{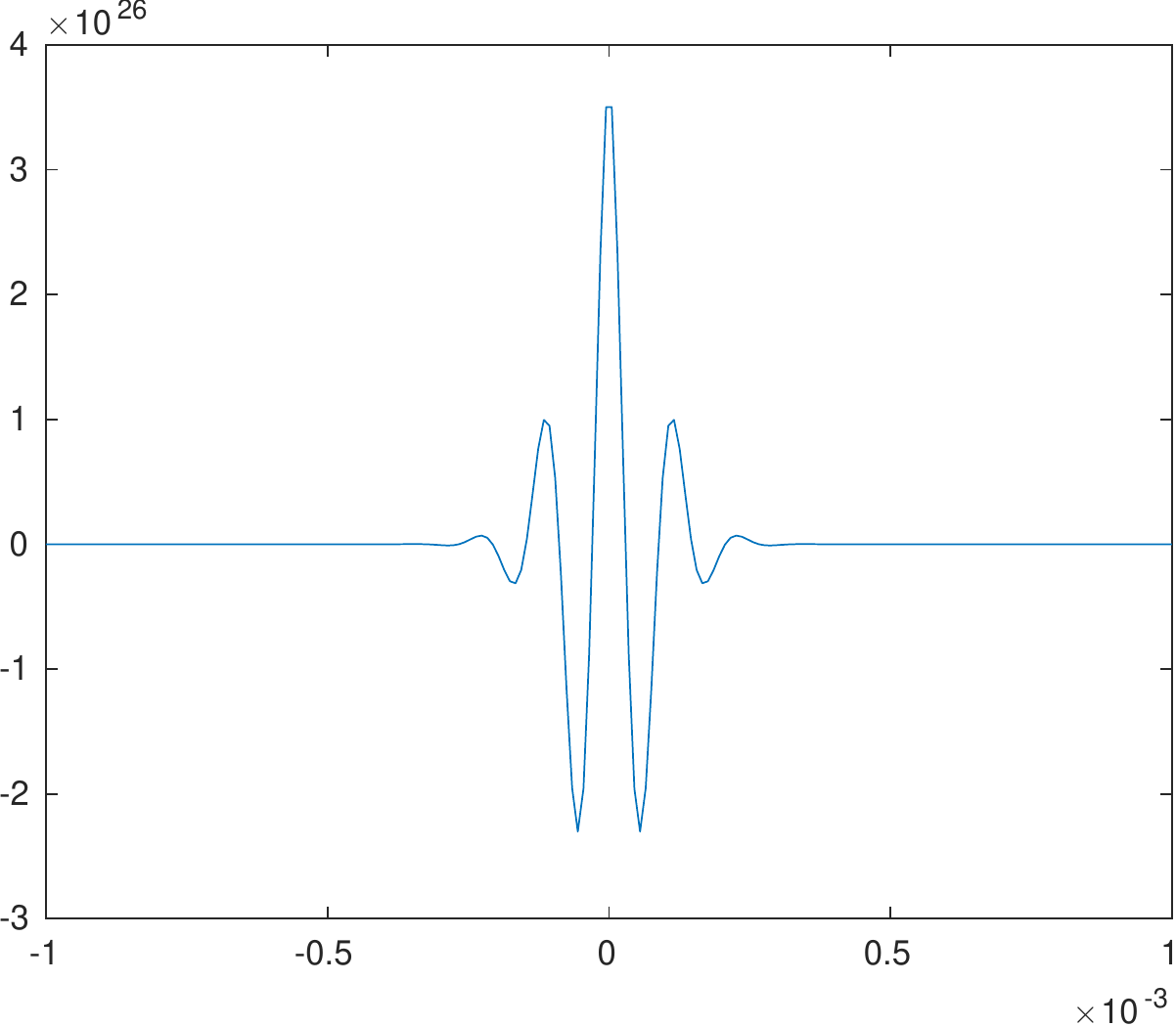}}
\vspace{2mm}
\subfloat[{Absolute values of the correlations in the clutter model ($\sigma=0$,
$v_c=\unit[10^{-2}]{m\!\cdot\! s^{-1}}$) and in the blood model ($\sigma^{2}=\unit[10^{-6}]{m^{2}s^{-1}}$,
$v_b=\unit[10^{-2}]{m\!\cdot\! s^{-1}}$).}\label{fig:corr_stat}]{\includegraphics[width=.98\textwidth]{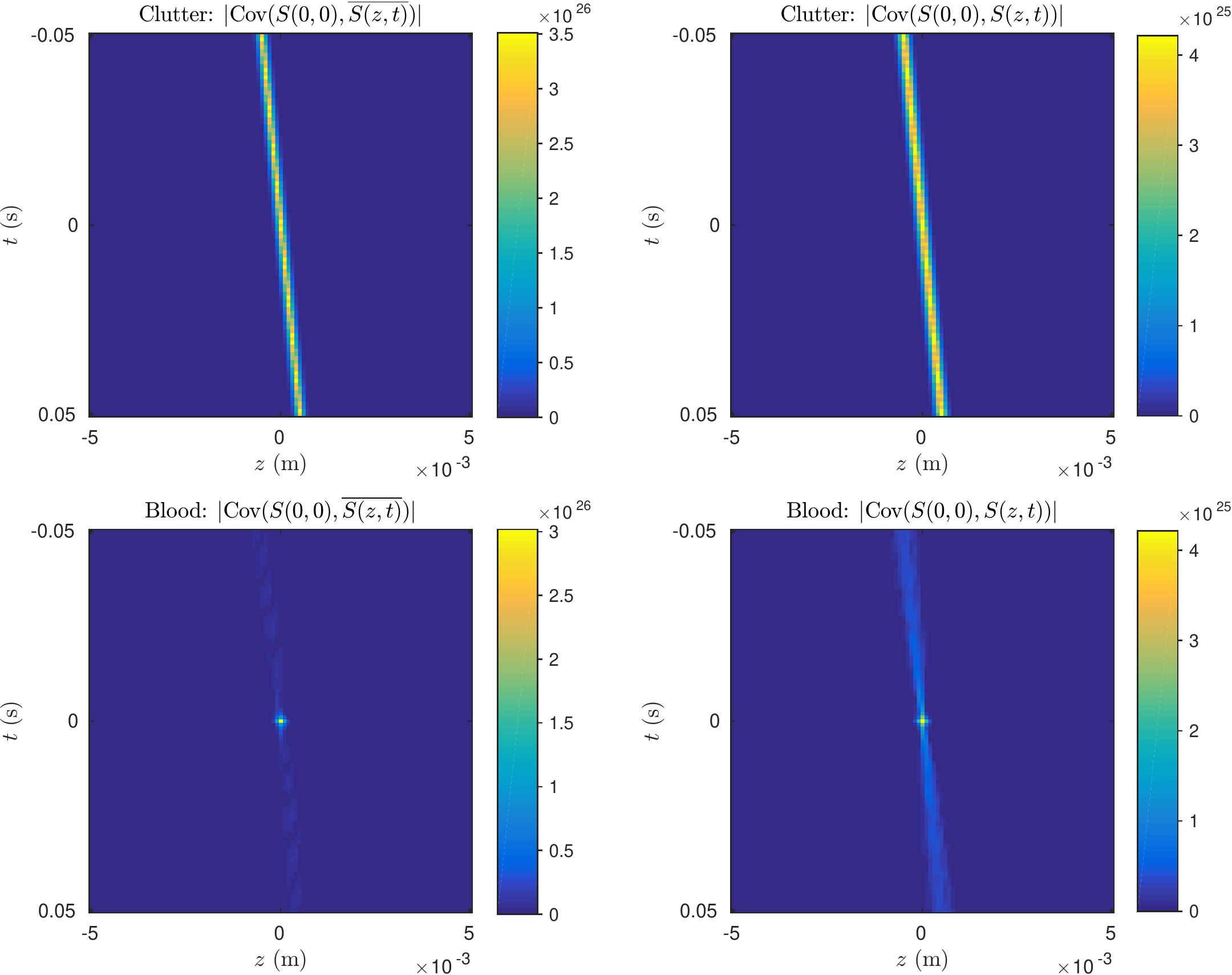}}

\caption{Correlations of the Casorati matrix.}
\end{figure}
On one hand, in the case of blood, since $C_{g\bar{g}}$ and $C_{gg}$
are oscillating and with very small support (see Figures~\ref{fig:cgg} and \ref{fig:cgg2}), the integration done
when taking the expectation in the blood case should yield
small correlations as long as $|t_{j'}-t_{j}|$ is large enough. On the
other hand, in the case of clutter, correlations will be high between
the two signals as long as $z_{i}-z_{i'}$ and $v_{c}\left(t_{j}-t_{j'}\right)$
are of the same order and almost cancel out. This heuristic is confirmed
by numerical experiments. In Figure~\ref{fig:corr_stat}, we compare
the clutter model and the blood model in one dimension: velocities
are in the $z$ direction, and we only consider points aligned on
the $z$ axis. As we can see, correlations are quickly decaying as
we move away from $(0,0)$ in the case of blood. In the case of clutter,
there are correlations at any times at the corresponding displaced locations.

\begin{figure}
\captionsetup[subfigure]{width=170pt}
\subfloat[{\label{fig:dist_eigs}The clutter model
($\sigma=0$, $v_c=\unit[10^{-2}]{m\!\cdot\! s^{-1}}$), the blood model
($\sigma^{2}=\unit[10^{-6}]{m^{2}s^{-1}}$, $v_b=\unit[10^{-2}]{m\!\cdot\! s^{-1}}$) and a white noise model with same variance as the blood.}]{\begin{centering}
\includegraphics[width=0.47\textwidth]{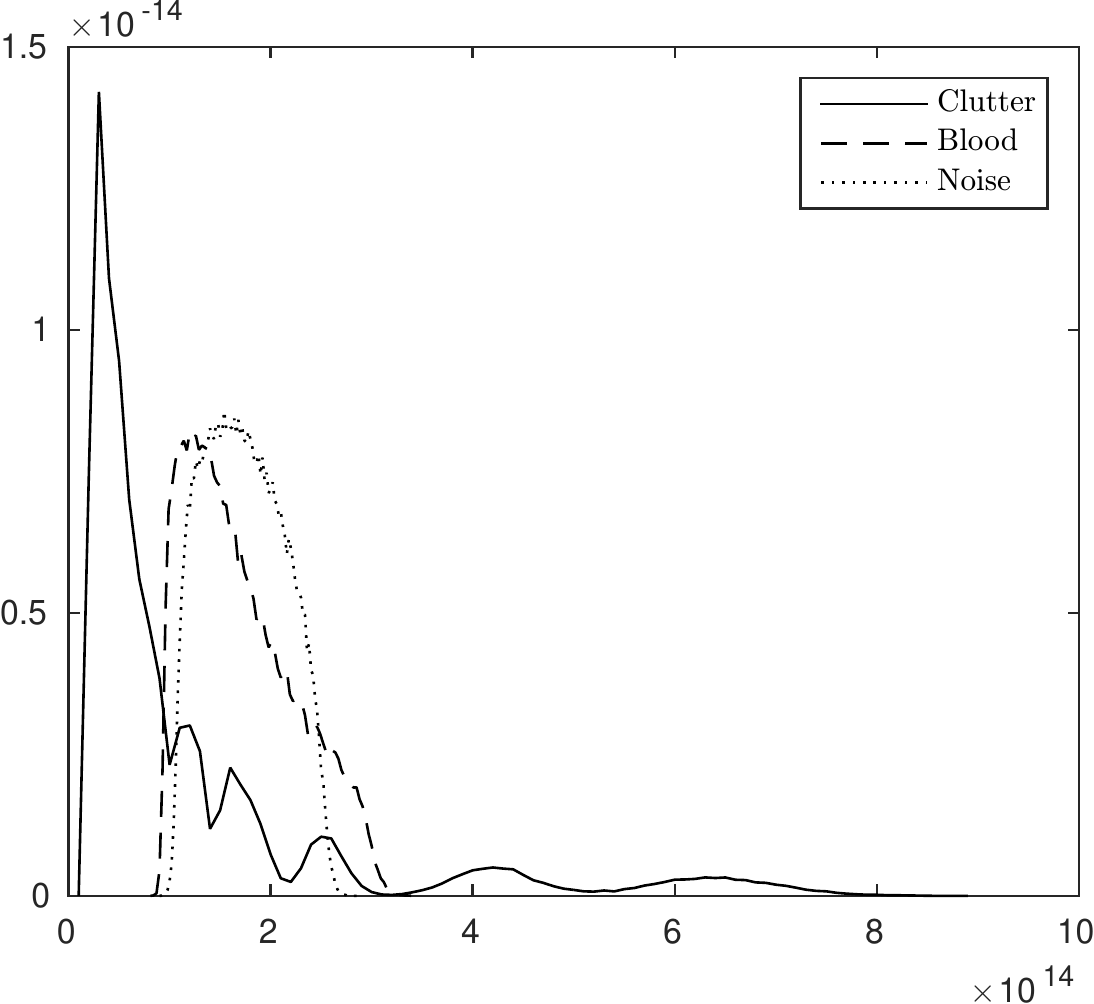}
\par\end{centering}

}\subfloat[\label{fig:dist_eigs-1}The clutter
model with different velocities.]{\begin{centering}
\includegraphics[width=0.47\textwidth]{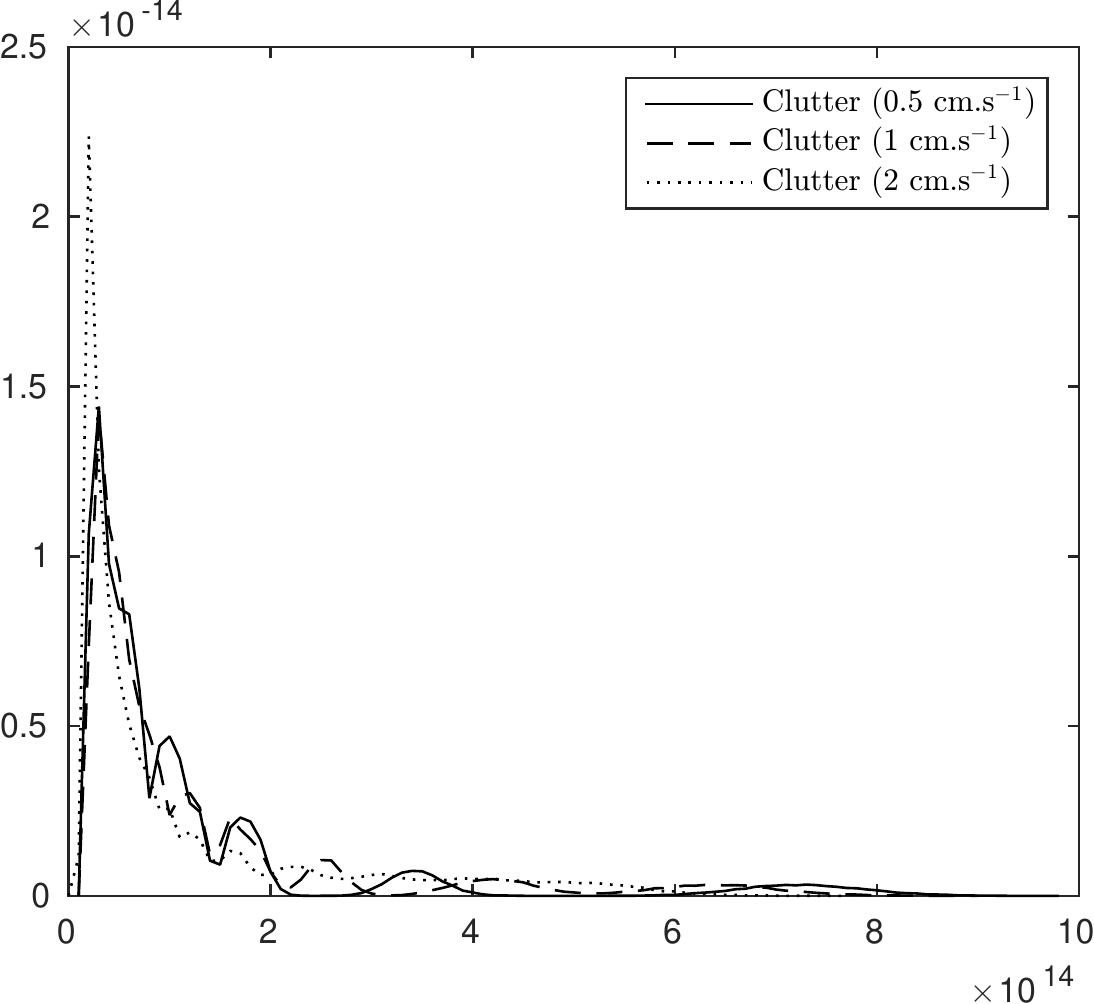}
\par\end{centering}

}\caption{The distribution of the singular values of the Casorati matrix $S$ in different cases.}

\end{figure}

Once the correlation matrix is computed, we can generate a large number
of samples to study the distribution
of the singular values in different cases. In Figure~\ref{fig:dist_eigs},
we compare the distribution in the two models (blood and clutter), using the 
Gaussian limit approximation for the simulations, with the same intensity for both models. A comparison with a
white noise model with the same variance shows that blood and noise have
approximately the same singular value distribution. On the contrary,
the distribution of the singular values of clutter presents a much larger
tail. A comparison of the distribution of the singular values for the clutter
model at different velocities shows no real difference in the tail
of the distribution  (Figure~\ref{fig:dist_eigs-1}).

As a consequence, the clutter signal $s_{c}$ is well
approximated by a low rank matrix, and the blood signal can be thought
of as if it were only noise. Therefore, the SVD method act as a denoising algorithm and extracts
the clutter signal, according to the discussion in the previous subsection.

\section{Numerical Experiments}  \label{sect5}

In this section, we consider again a more realistic 2D model, given by \eqref{eq:a_k}. This framework
will allow us to simulate generic blood flow imaging sequences from
particles. The dynamics of blood and clutter are modelled as follows.
Let us assume that clutter is subject to a deterministic and computable
flow $\varphi_{c}$. The randomness of the motion of red blood cells
in vessels is modelled by a stochastic differential equation, given by
\begin{equation}
d\mathbf{y}=\mathbf{v}_{b}\left(t,y\right)dt+\sigma (y)\,dB_t,\label{eq:sde}
\end{equation}
where $B_t$ is a two dimensional Brownian motion and $\sigma$
is determined by the effective diffusion coefficient $K=\frac{1}{2}\sigma^{2}$.
In blood vessels, this diffusion coefficient is proportional to the
product $\dot{\gamma}r^{2}$ where $\dot{\gamma}$ is the shear stress
in the vessel, and $r$ is the radius of red blood cells. As in the
previous section, let $a_{c}=a_{1,c}$ and $a_{b}=a_{1,b}$. Let $\varphi_{b}$
be the flow associated to \eqref{eq:sde}. We assume that
$\varphi_{b}$ represents the dynamics of blood particles, relative
to overall clutter movement, so that
\begin{equation}
a_{c}\left(t\right)=\varphi_{c}\left(u_c,t\right),\qquad \varphi_c(u_c,0)=u_c,\label{eq:clutdyn}
\end{equation}
and
\begin{equation}
a_{b}\left(t\right)=\varphi_{c}\left(\varphi_{b}\left(u_b,t\right),t\right),\qquad \varphi_b(u_b,0)=u_b.
\label{eq:blooddyn}
\end{equation}
The dynamics of all the other particles are then taken to be independent
realizations of the same dynamics. The velocity field $\mathbf{v}_{b}$
and the clutter dynamics $\varphi_{c}$ are computed beforehand and
correspond to the general blood flow velocity and to an elastic displacement, respectively. In our experiments, we let $\varphi_{c}$ be
an affine displacement of the medium, changing over time: a global affine transformation, with slowly varying translation and shearing 
applied to the medium at each frame, namely
\[
\varphi_c(u,t)=\left[\begin{smallmatrix} 1 & w_1(t) \\ 0 & 1 \end{smallmatrix}\right]u + \left[\begin{smallmatrix} w_2(t) \\ w_3(t) \end{smallmatrix}\right],
\]
where $w_i$ are smooth and slowly varying (compared to $\varphi_b$) functions such that $w_i(0)=0$. As
for the blood velocity flow $\mathbf{v}_{b}$, it is parallel to the blood vessels,
with its intensity decreasing away from the center of the blood vessel \cite[Section~11.3]{szabo_diagnostic_2013}. More precisely,
 $\mathbf{v}_{b}$ is a Poiseuille laminar flow, namely the mean blood flow velocity is half
of the maximum velocity, which is the fluid velocity in the center of the vessel. 

\begin{figure}
\begin{centering}
\includegraphics[width=0.4\textwidth]{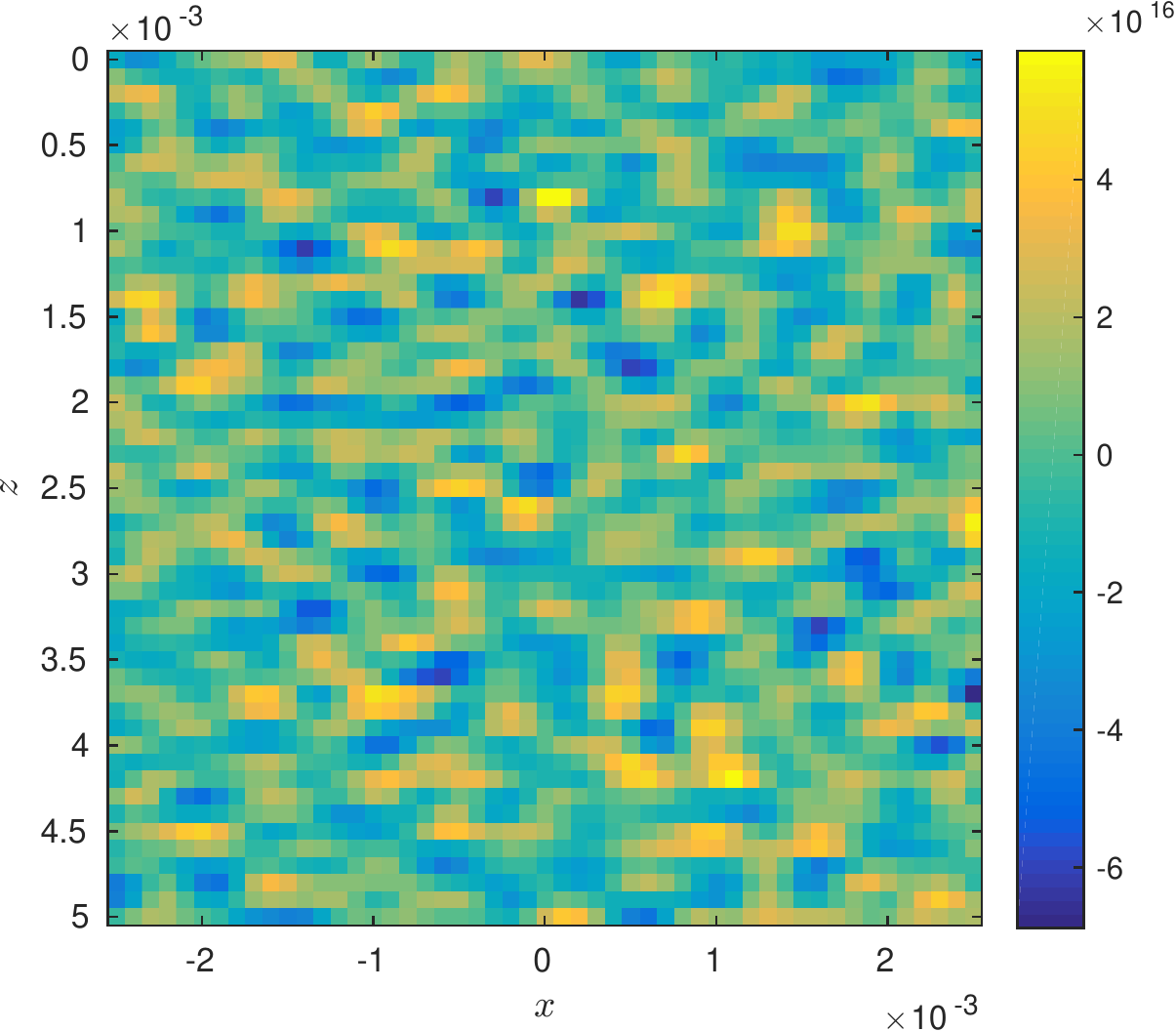}
\par\end{centering}
\caption{\label{fig:singleframe}Single frame of ultrafast ultrasound (real part).}
\end{figure}

The relative blood displacements $b_{k,j}=\varphi_{k,b}\left(u_{b,k},t_j\right)$
are computed according to the following discretization of the stochastic
differential equation~\eqref{eq:sde}:
\[
b_{k,j+1}=b_{k,j}+\delta t\mathbf{v}_{b}\left(t_{j},b_{k,j}\right)+\sqrt{\delta t}\sigma\left(b_{k,j}\right)X_{k,j}+o\left(\delta t\right),
\]
where $\left(X_{k,j}\right)$ are centered independent Gaussian random
variables and $\delta t=t_{j+1}-t_j$ is taken to be constant. The blood particle positions $a_{k,b}\left(t_{j}\right)$
are then computed simply by applying the precomputed flow $\varphi_{c}$.

\begin{figure}
\centering
\subfloat[{Maximum blood velocity: $\unit[2]{cm\!\cdot\! s^{-1}}$; mean clutter velocity: $\unit[1]{cm\!\cdot\! s^{-1}}$.}]
{ \includegraphics[width=0.84\textwidth]{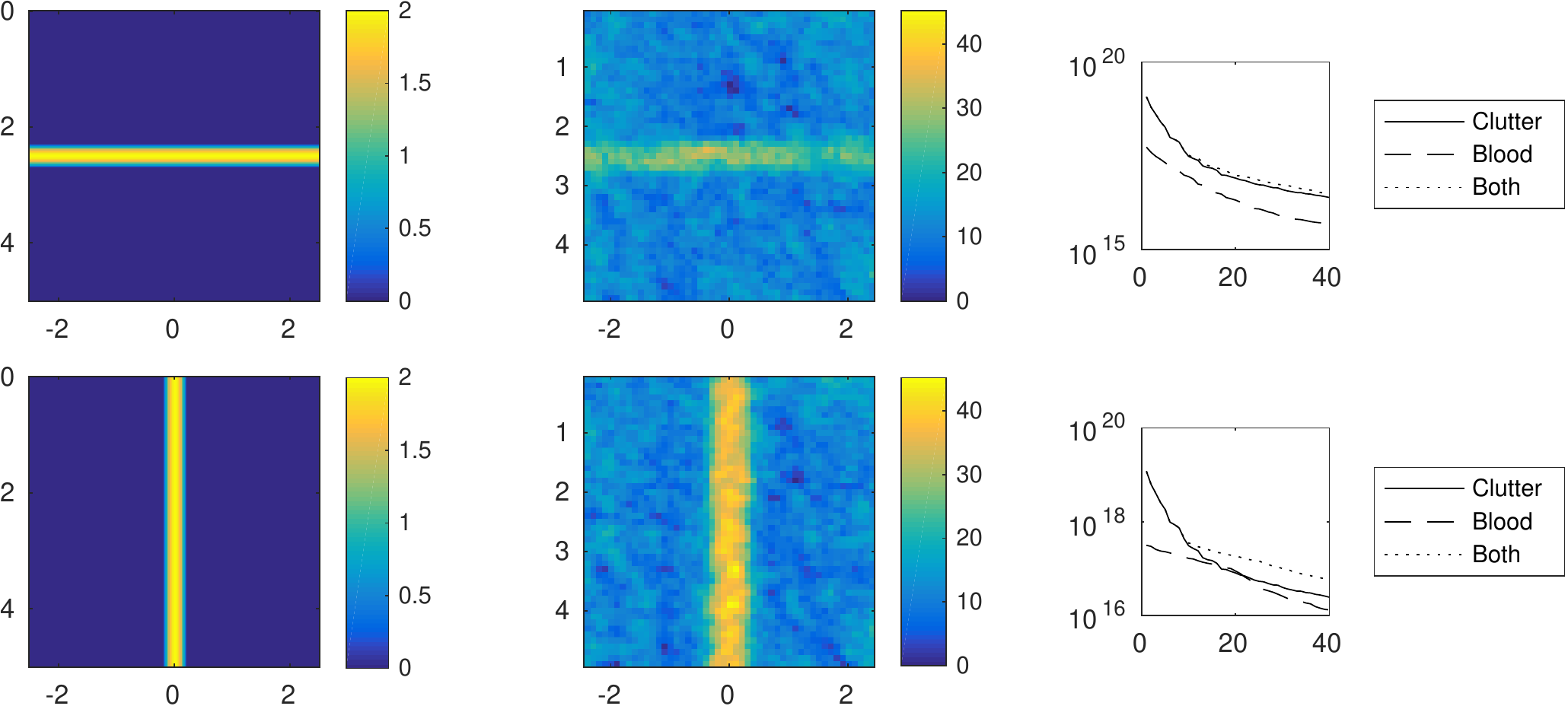} }

\subfloat[{Maximum blood velocity: $\unit[1]{cm\!\cdot\! s^{-1}}$; mean clutter velocity: $\unit[1]{cm\!\cdot\! s^{-1}}$.}]
{
\includegraphics[width=0.84\textwidth]{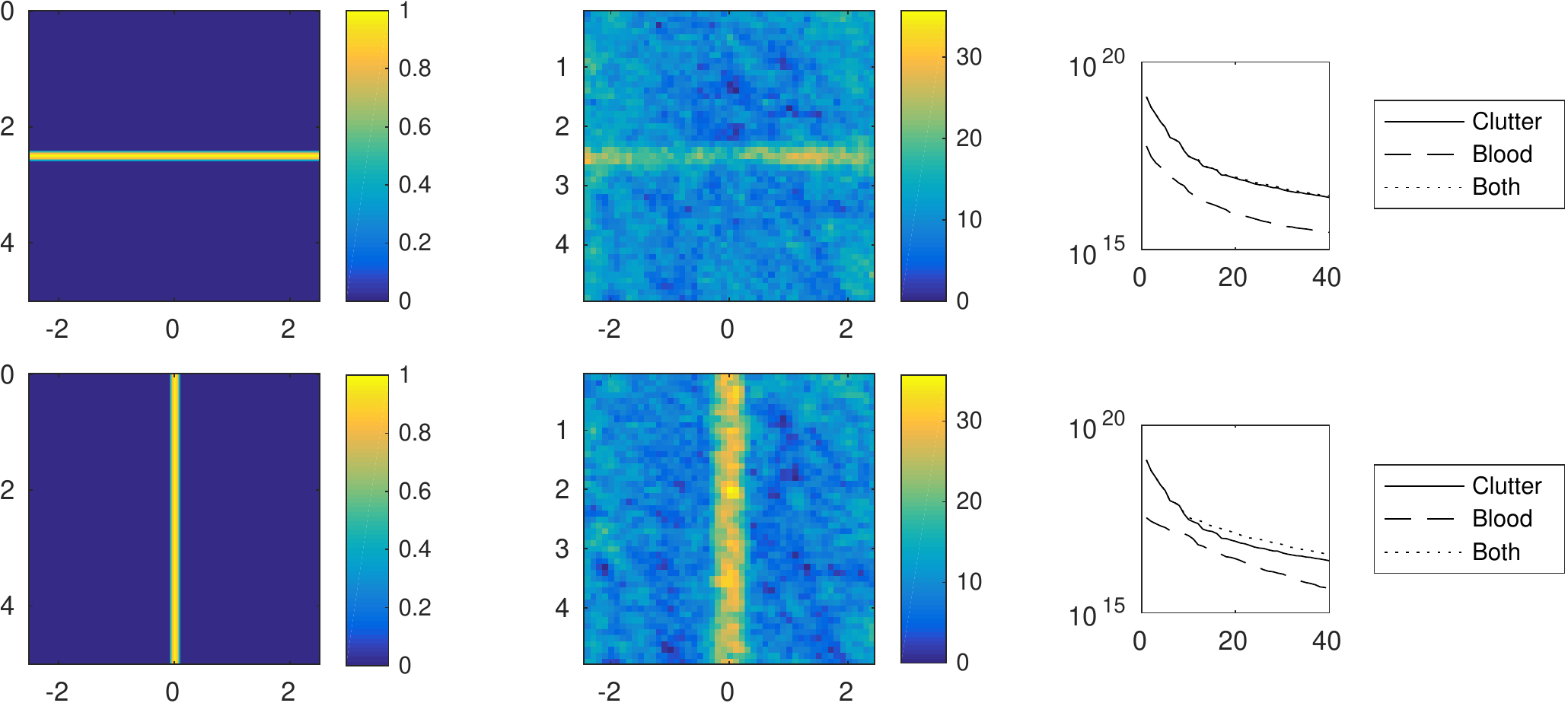}
\label{fig:SVD-2d-b}
}

\subfloat[{Maximum blood velocity: $\unit[0.5]{cm\!\cdot\! s^{-1}}$; mean clutter velocity: $\unit[1]{cm\!\cdot\! s^{-1}}$.}]
{
\includegraphics[width=0.84\textwidth]{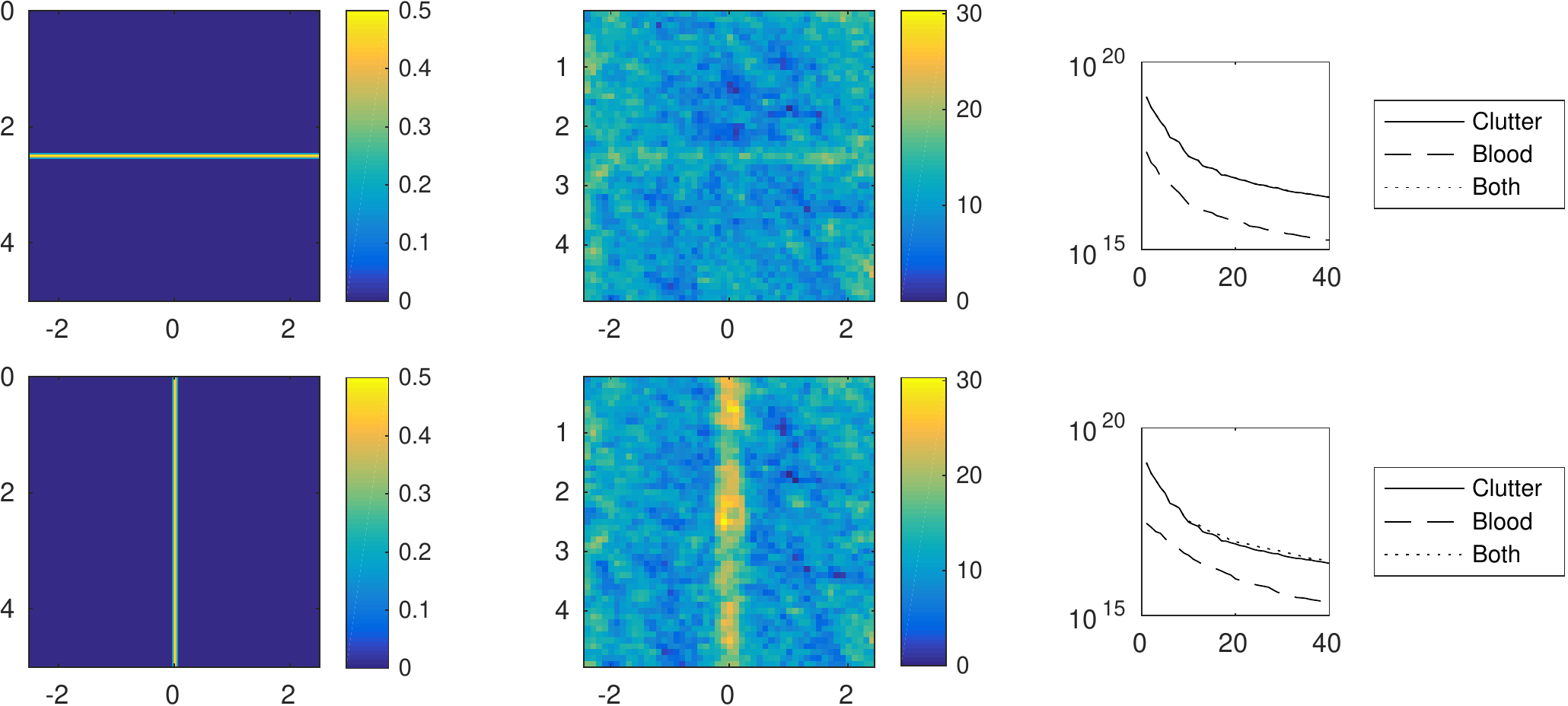}
\label{fig:SVD-2d-c}
}

\caption{The SVD method for different velocities and orientations. In each case,
we have from left to right: the
blood velocity and location, the reconstructed blood location, the decay of the singular values. The squares  are $\unit[5]{mm}\times\unit[5]{mm}$, and the horizontal and vertical axes are the $x$ and $z$ axes, respectively. The parameters used are those given in \eqref{eq:f-chi} and  \eqref{eq:quantities}, $F=0.4$ and $\Theta=7°$. The density of particles for both blood and clutter is 2,000 per $\unit{mm^2}$, and $\sigma = 2.5\cdot 10^{-5}$.\label{fig:SVD-2d}} 
\end{figure}

In order to validate the SVD approach, we explore the effects  of the blood
velocity and of the direction of the blood vessels on the behavior of the singular
values and on the quality of the reconstruction. In each case, the clutter displacement is the same composition of time-varying
shearing and translation, and the mean clutter velocity is $\unit[1]{cm\!\cdot\! s^{-1}}$. We choose $C_c = 5$ and $C_b = 1$, for the same density of scatterers
from clutter and blood: per unit of area, the clutter intensity  is therefore five times higher than the blood intensity. A single frame of ultrafast
ultrasound imaging is presented in Figure~\ref{fig:singleframe}: it is clear that without further processing, it is impossible to locate the blood vessels.

\begin{figure}
\hspace{7pt}\subfloat[Flow parallel to the receptor array.]{\includegraphics[width=0.45\textwidth]{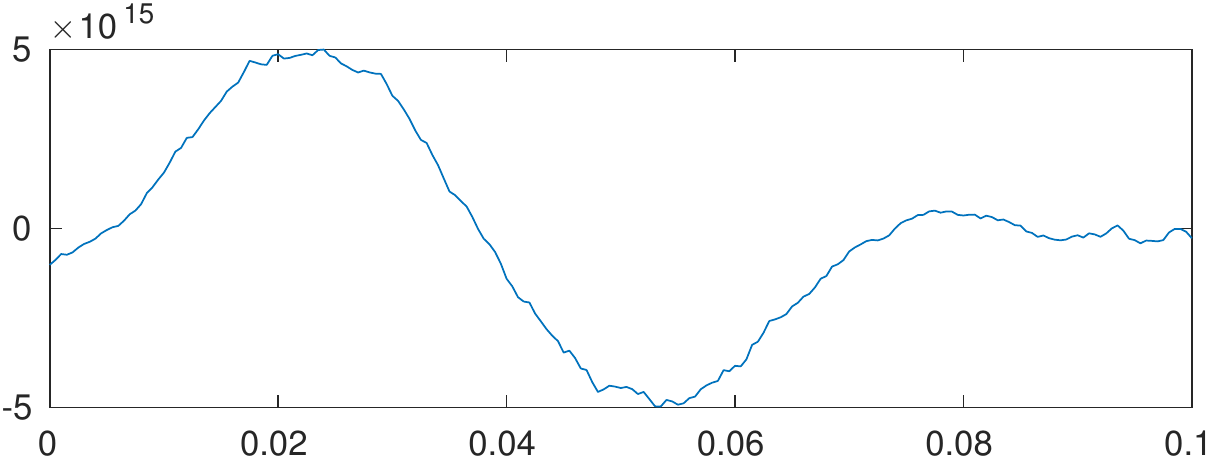}
}\hspace{13pt}\subfloat[Flow perpendicular to the receptor array.]{\includegraphics[width=0.45\textwidth]{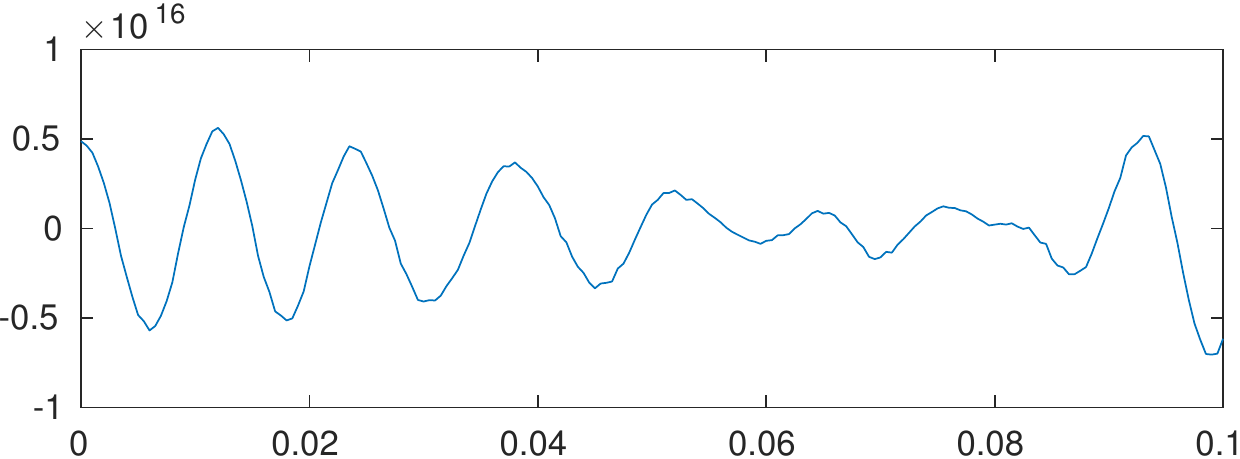}
}
\caption{\label{fig:singlepixel}Time behavior of a single pixel (real part),
located in a constant velocity flow.}
\end{figure}

In Figure~\ref{fig:SVD-2d},
the results for various velocities and orientations are presented.
The reconstruction intensities
are expressed in decibels, relatively to the smallest value in the image.
The SVD method allows for reconstruction
of blood vessels, even if the maximum blood velocity is close to, or oven lower than, the mean velocity of clutter.
We always use the threshold $K=20$.
As we can see, due to the better resolution in the $z$ direction discussed in Section~\ref{sect2},
vessels oriented parallel to the receptor array have a reconstruction
with a better resolution. But due to the oscillating behavior of
the PSF in the $z$ direction, and the low-pass filter behavior of
the PSF in the $x$ direction, the sensitivity is better for vessels
oriented perpendicularly to the receptor array, and the SVD method is
able to reconstruct smaller vessels with lower velocities. This follows from the discussion in Subsection~\ref{sub:doppler}. In order to visualize this phenomenon even better,  Figure~\ref{fig:singlepixel} presents the time behavior of
a single pixel from the data of Figure~\ref{fig:SVD-2d-c}. We can clearly see the doppler effect in the case
when the flow is perpendicular to the receptor array, and the low
frequency behavior of the signal in the case when it is parallel
to the receptor array.

\begin{figure}
\centering
\subfloat[Flow parallel to the receptor array.]{\includegraphics[width=\textwidth]{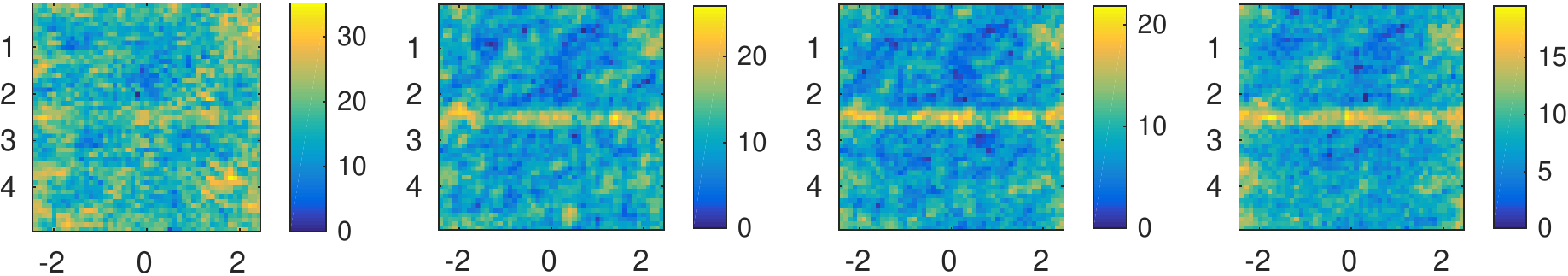}}

\subfloat[Flow perpendicular to the receptor array.]{\includegraphics[width=\textwidth]{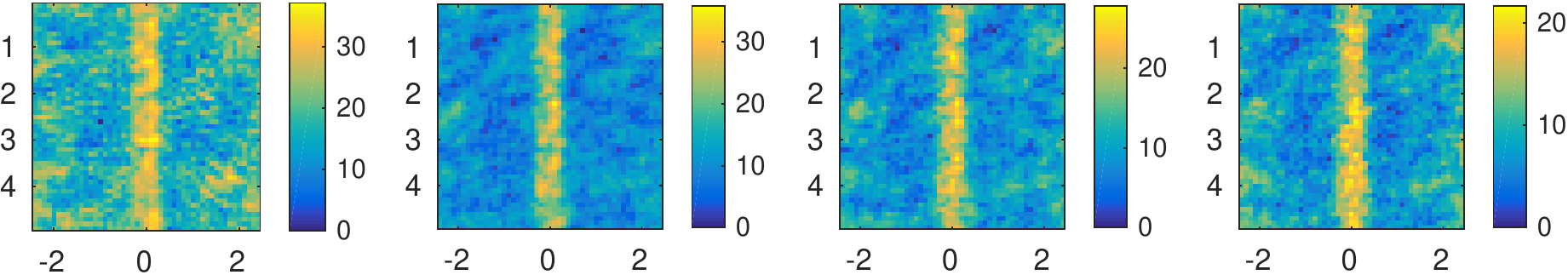}}
\caption{\label{fig:threshold}Effect of the threshold $K$ on the reconstruction. From left to right: ${K=10,20,30,40}$.}
\end{figure}

In Figure~\ref{fig:threshold}, results of an investigation on the effect of the threshold $K$ on the reconstruction are presented. Except for 
$K$, the parameters of Figure~\ref{fig:SVD-2d-b} are used. If the threshold is too low,
the reconstruction is not satisfactory and artefacts appear everywhere in the reconstructed image. If the threshold is too high, the reconstruction still works
but the contrast becomes lower. With our parameters, $K=20$ seems to produce the best results.

\begin{figure}
\centering
\subfloat[Contrast as a function of noise level]{\includegraphics[width=.68\textwidth]{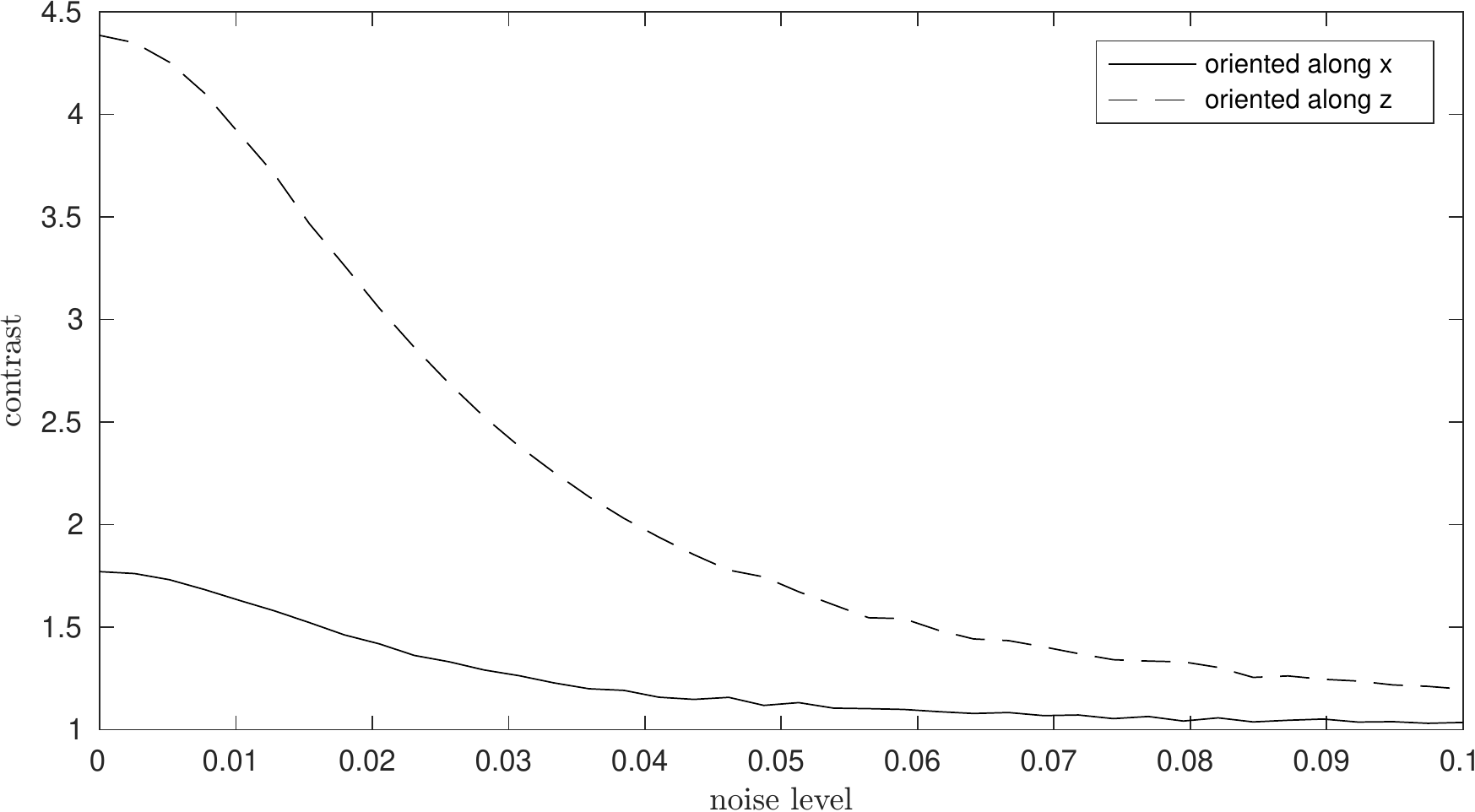}
}

\subfloat[0\% noise]{\includegraphics[width=0.235\textwidth]{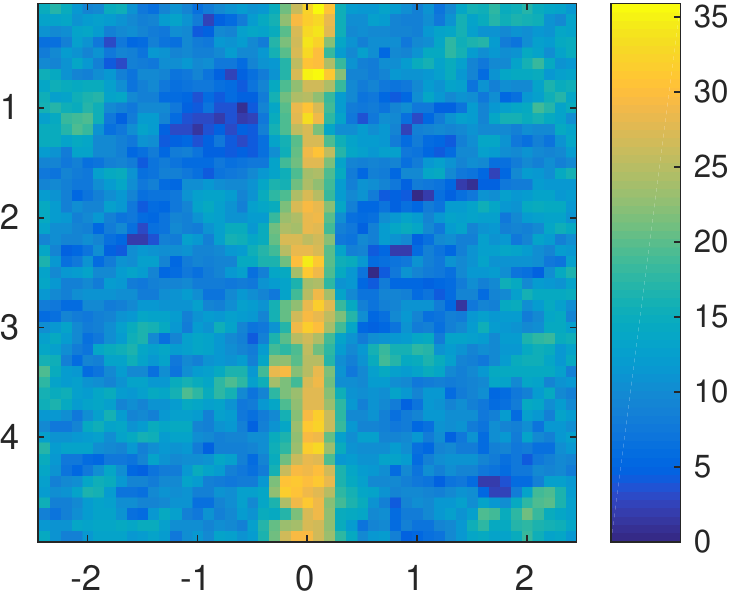}
}\hfill\subfloat[2.5\% noise]{\includegraphics[width=0.235\textwidth]{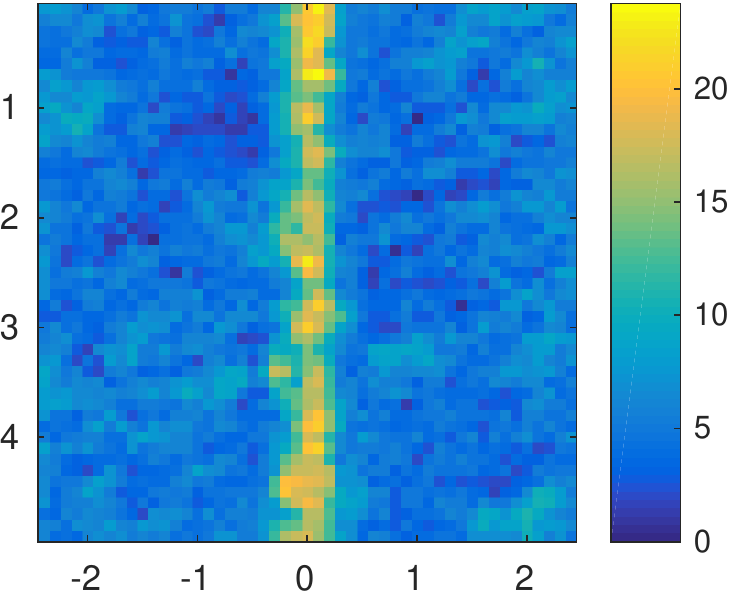}
}\hfill\subfloat[5\% noise]{\includegraphics[width=0.235\textwidth]{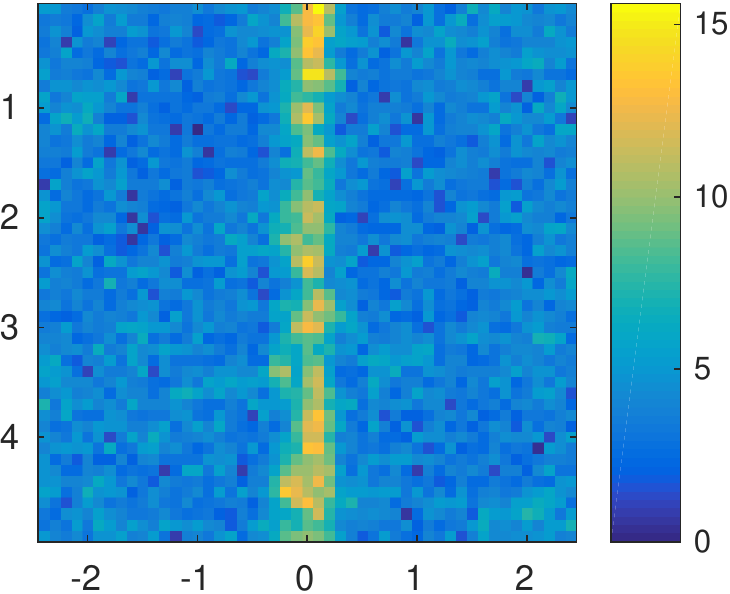}
}\hfill\subfloat[7.5\% noise]{\includegraphics[width=0.235\textwidth]{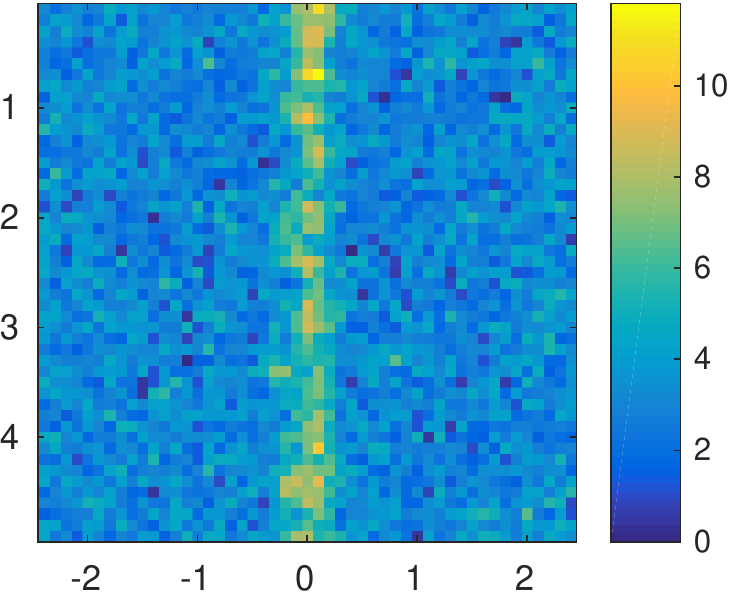}
}
\caption{\label{fig:contrast_noise}Effect of noise on the reconstruction. The parameters are the same used in Figure~\ref{fig:SVD-2d-b}.}
\end{figure}

In order to further validate the method, we consider the impact of measurement noise on the recovery. To this end, we add independent 
white Gaussian noise to the data, and consider the quality of the reconstruction as a function of the noise intensity.
Let us define the contrast of the reconstruction as the ratio between the mean intensity of the reconstructed image 
inside and outside the blood domain. The parameters of Figure~\ref{fig:SVD-2d-b} are used. Blood 
intensity is five times lower than clutter intensity, and therefore a noise intensity of 10\% corresponds to half the intensity 
of blood. In Figure~\ref{fig:contrast_noise}, sample reconstructions at different noise levels are provided. We can conclude that contrast is robust to moderate 
levels of noise, since blood vessels can still be identified up to 7.5\% of noise if they are oriented along the z axis, and up to
2.5\% of noise if they are oriented along the x axis. Figure~\ref{fig:contrast_noise} also clearly quantifies the better contrast
for vessels oriented along the z axis.

\section{Concluding Remarks}  \label{sect6}

In this paper, we have provided for the first time a detailed mathematical analysis of ultrafast ultrasound  imaging. By using a random model for the movement of the blood cells, we have shown  that an SVD approach can separate the blood signal  from the clutter signal. Our model and results  open a door for a mathematical
and numerical framework for realizing super-resolution in dynamic optical coherence tomography \cite{dynamicoct},
in ultrafast ultrasound imaging by tracking microbubbles \cite{Errico2015},
as well as in acousto-optic imaging based on the use
of ultrasound plane waves instead of focused ones, which allows to increase
 the imaging rate drastically \cite{acousto-optic}. These three modalities are under investigation and their mathematical and numerical modeling will be the subject of forthcoming papers.

\appendix

\section{\label{sec:The-justification-of}The Justification of the Approximation
of the PSF}

This appendix is devoted to the formal justification of the PSF approximation
(\ref{eq:g-conv}) which was obtained by truncating the Taylor expansion
of $w_{\pm}^{\theta}$ at the first order: we shall show here that
the error caused by this truncation is small. For simplicity, we shall
consider only the case when $z=z'$ and $\theta=0$: the general case
may be tackled in a similar way. Without loss of generality, we may
set $x'=0$ and suppose $x\ge0$. We also suppose that we are not
too close to the detectors, namely $z\ge\unit[10^{-2}]{m}$. Moreover,
in order to be able to be quantitative, we consider the particular
case when $F=0.4$ and $\tau=1$.

The expression of the PSF that we want to approximate is (see (\ref{eq:g-app1}))
\[
g(x):=g_{0}((x,z),(0,z))=\frac{c_{0}}{4\pi x}\left[f'(w_{+}(x))-f'(w_{-}(x))\right],
\]
where $w_{\pm}(x)$ is given by
\[
w_{\pm}(x):=h_{\x,\x'}^{0}(x\pm Fz)=c_{0}^{-1}\left(\sqrt{1+F^{2}}z-\sqrt{z^{2}+(x\pm Fz)^{2}}\right).
\]
(Note that, for simplicity of notation, we have removed the dependence
of $w$ on $\theta$ and $z$.) An immediate calculation shows that
\[
w_{\pm}(0)=0,\qquad w'_{\pm}(0)=\frac{\mp c_{0}^{-1}F}{\sqrt{1+F^{2}}},\qquad w''_{\pm}(x)=\frac{-c_{0}^{-1}z^{2}}{\left((x\pm Fz)^{2}+z^{2}\right)^{3/2}}.
\]
Hence, there exists $\xi_{x}\in[0,x]$ such that
\[
w_{\pm}(x)=\frac{\mp c_{0}^{-1}F}{\sqrt{1+F^{2}}}x+c_{x}\frac{x^{2}}{2},\qquad|c_{x}|=|w_{\pm}''(\xi_{x})|\le c_{0}^{-1}z^{-1}.
\]
Therefore, the absolute error $E(x)$ due to the truncation of the
Taylor series of $w_{\pm}$ at first order is given by
\[
E(x)=c_{0}(4\pi)^{-1}\left[E_{+}(x)-E_{-}(x)\right],
\]
where
\[
E_{\pm}(x)=\frac{1}{x}\left[f'(\frac{\mp c_{0}^{-1}F}{\sqrt{1+F^{2}}}x+c_{x}\frac{x^{2}}{2})-f'(\frac{\mp c_{0}^{-1}F}{\sqrt{1+F^{2}}}x)\right].
\]

We now consider two cases, depending on $x$. First, consider the
case when $x>\unit[5\cdot10^{-3}]{m}$. From the above calculations
we immediately have
\[
|E(x)|\le c_{0}(4\pi)^{-1}\frac{4}{x}\left\Vert f'\right\Vert _{\infty}\le\frac{2}{5}c_{0}10^{3}\nu_{0}\le3.7\cdot10^{12}.
\]

Next, consider the case when $x\le\unit[5\cdot10^{-3}]{m}$. By using
again the mean value theorem we obtain
\[
E_{\pm}(x)=c_{x}\frac{x}{2}f''(\theta_{x}),\qquad\theta_{x}=\frac{\mp c_{0}^{-1}F}{\sqrt{1+F^{2}}}x+\delta_{x}c_{x}\frac{x^{2}}{2}
\]
for some $\delta_{x}\in[0,1]$. Since $|f''(t)|$ is even and decreasing
for $t>0$, we have that
\[
|E_{\pm}(x)|\le c_{0}^{-1}\frac{x}{2z}|f''(\frac{c_{0}^{-1}F}{\sqrt{1+F^{2}}}x-c_{0}^{-1}\frac{x^{2}}{2z})|,
\]
since the inequality $x\le\unit[5\cdot10^{-3}]{m}$ guarantees that
$\frac{c_{0}^{-1}F}{\sqrt{1+F^{2}}}x-c_{0}^{-1}\frac{x^{2}}{2z}>0$.
Therefore we have
\[
|E(x)|\le(4\pi)^{-1}xz^{-1}|f''(\frac{c_{0}^{-1}F}{\sqrt{1+F^{2}}}x-c_{0}^{-1}\frac{x^{2}}{2z})|.
\]
Let us look at the right hand side of this inequality. As $x\to0$
the error tends to $0$: this is expected, because of the Taylor expansion
around $0$. On the other hand, for big $x$, the value of $|f''(\frac{c_{0}^{-1}F}{\sqrt{1+F^{2}}}x-c_{0}^{-1}\frac{x^{2}}{2z})|$
is very small, since $|f''(t)|$ decays very rapidly for large $t$.
Therefore, the maximum of the right hand side is attained in a point
$x^{*}\in(0,0.005)$. The value in this point may be explicitly calculated,
and we have
\[
|E(x)|\le4\cdot10^{12},\qquad0\le x\le\unit[5\cdot10^{-3}]{m}.
\]

To summarize the above derivation, we have shown that the absolute
error $E(x)$ is bounded by
\begin{equation}
|E(x)|\le4\cdot10^{12},\qquad x\ge0.\label{eq:relative-error}
\end{equation}

We now wish to estimate the relative error $\left\Vert E\right\Vert _{\infty}/\left\Vert g\right\Vert _{\infty}$.
In order to do this, let us compute $g(0)$. Since the Taylor expansion
becomes exact as $x\to0$, we may very well compute $g(0)$ by using
the approximated version. Thus, setting $G=F/\sqrt{1+F^{2}}$ we have
\begin{equation*}
\begin{split}g(0) & =\lim_{x\to0}-\frac{c_{0}}{4\pi x}\left[f'(c_{0}^{-1}Gx)-f'(-c_{0}^{-1}Gx)\right]\\
 & =\lim_{x\to0} -G(4\pi)^{-1}\left[\frac{f'(c_{0}^{-1}Gx)-f'(0)}{c_{0}^{-1}Gx}+\frac{f'(-c_{0}^{-1}Gx)-f'(0)}{-c_{0}^{-1}Gx}\right]\\
 & =-2G(4\pi)^{-1}f''(0),
\end{split} \end{equation*}
whence $|g(0)|\ge8.8\cdot10^{13}$ by a direct calculation of $|f''(0)|$.
Finally, combining this inequality with (\ref{eq:relative-error})
allows to bound the relative error by
\[
\frac{\left\Vert E\right\Vert _{\infty}}{\left\Vert g\right\Vert _{\infty}}\le5\%.
\]

We have proven that the relative error of the approximation obtained
by truncating the Taylor expansions of $w_{\pm}$ at the first order
is less than 5\%. This has been proven only in the particular case
when $z=z'$: the general case may be done by extending the above
argument to two dimensions.

\bibliographystyle{abbrv}
\bibliography{ultrafast}

\end{document}